\documentclass[a4paper,10pt]{amsart}

\usepackage{amsmath}
\usepackage{amsthm}
\usepackage{amssymb}
\usepackage{amsfonts}
\usepackage{mathrsfs}  
\usepackage{enumitem} 
\usepackage{xcolor}
\usepackage[bookmarks=false,hyperindex,pdftex,colorlinks,citecolor=blue,urlcolor=cyan]{hyperref}
\usepackage[a4paper,lmargin=3cm,rmargin=3cm,tmargin=4cm,bmargin=4cm,marginparwidth=2.8cm,marginparsep=1mm]{geometry} 
\usepackage{soul} 
\usepackage{marginnote} 
\usepackage{tikz}

\DeclareMathOperator{\diam}{diam}                           
\DeclareMathOperator{\Lip}{Lip}                             
\DeclareMathOperator{\SEP}{SEP}
\DeclareMathOperator{\STD}{STD}

\newcommand{\NN}{\mathbb{N}}                                
\newcommand{\ZZ}{\mathbb{Z}}                                
\newcommand{\RR}{\mathbb{R}}                                
\newcommand{\TT}{\mathbb{T}}                                
\newcommand{\EE}{\mathbb{E}}

\newcommand{\abs}[1]{\left|{#1}\right|}                     
\newcommand{\pare}[1]{\left({#1}\right)}                    
\newcommand{\set}[1]{\left\{{#1}\right\}}                   
\newcommand{\norm}[1]{\left\|{#1}\right\|}                  
\newcommand{\restrict}{\mathord{\upharpoonright}}           

\newcommand{\lipfree}[1]{\mathcal{F}({#1})}                 

\renewcommand{\leq}{\leqslant}
\renewcommand{\geq}{\geqslant}

\theoremstyle{plain}
\newtheorem{theorem}{Theorem}[section]
\newtheorem{lemma}[theorem]{Lemma}
\newtheorem{corollary}[theorem]{Corollary}
\newtheorem{proposition}[theorem]{Proposition}

\newtheorem*{claim*}{Claim}

\theoremstyle{definition}
\newtheorem*{definition*}{Definition}
\newtheorem{definition}[theorem]{Definition}

\theoremstyle{remark}

\begin{document}

\title{Nagata dimension, separation modulus and stochastic embeddings into trees}

\author[R. J. Aliaga]{Ram\'on J. Aliaga}
\address[R. J. Aliaga]{Instituto Universitario de Matem\'atica Pura y Aplicada,
Universitat Polit\`ecnica de Val\`encia,
Camino de Vera S/N,
46022 Valencia, Spain}
\email{ramon.aliaga@upv.es}

\author[R. Medina]{Rub\'en Medina}
\address[R. Medina]{Universidad P\'ublica de Navarra (ETSIIIT). Departamento de Estad\'istica, Inform\'atica y Matem\'aticas. Institute for Advanced Materials and Mathematics. 31009-Navarra (Spain). \newline
\href{https://orcid.org/0000-0002-4925-0057}{ORCID: \texttt{0000-0002-4925-0057}}}
\email{ruben.medina@unavarra.es}

\begin{abstract}
The separation modulus of a metric space, introduced by Bartal in 1996, is related to different notions of dimensionality. Most generally, Naor and Silberman gave in 2011 an upper bound on the separation modulus in terms of Nagata dimension. We prove that this bound cannot be reversed; that is, there are metric spaces with arbitrarily high Nagata dimension and stochastically embeddable into trees with uniform distortion, which therefore have uniformly bounded separation modulus. Similarly, we construct spaces with arbitrarily high Nagata dimension whose spaces of probability measures endowed with the Earth Mover's Distance (equivalently, their Lipschitz-free spaces) are embeddable into $L_1$  with uniformly bounded distortion. By gluing these spaces together, we obtain infinite Nagata dimensional metric spaces whose separation modulus is finite and whose Lipschitz-free space is isomorphic to $\ell_1$, respectively. Finally, we provide relations between separation modulus and stochastic embeddability into trees, also known as random approximability by tree metrics, in particular for the case of snowflake metrics.
\end{abstract}



\maketitle

\section{Introduction}

In this paper, we study the relation between several concepts in the theory of metric spaces that somehow measure the (finite) dimensionality of a metric space, and obtain applications to stochastic embeddings into trees, also known as random approximation by tree metrics, as well as embedding of Lipschitz-free spaces into $L_1$. We consider two different classes of metric notions:

\begin{itemize}
    \item Random decomposition moduli, in particular the separation modulus (Definition \ref{def:sep}). It provides a measure of how likely it is to randomly find a `well-behaved' partition of a metric space, that is, a partition that separates a previously chosen pair of points, in comparison with the distance between them.

    \item Metric dimensions, in particular Nagata dimension (Definition \ref{def:nagata}) and Lipschitz dimension (Definition \ref{def:lipdim}): the former is a metric version of topological dimension in terms of covers that must be valid at all scales, and the latter is expressed in terms of existence of sufficiently nice maps into Euclidean space. Both are biLipschitz invariant.
\end{itemize}

Our work is motivated by the existence of several results in the literature involving Nagata dimension that we do not consider completely satisfactory for a number of reasons: for some of them Nagata dimension may not be the best possible notion to consider, and for others the validity of the converse remains open. We focus on four specific aspects and attempt to improve on their current understanding. We now list them all, together with our main results.

\subsection*{1. Upper bounds for the separation modulus}
In \cite[Theorem 1.3]{Naor17}, Naor shows that the separation modulus of a $d$-dimensional normed space $X_d$, $\SEP(X_d)$, is bounded above by a universal multiple of $d$. One may also find bounds for general metric spaces: by \cite[Lemma 5.1]{NS11} and \cite[Theorem 2.2]{LN25}, if a metric space $M$ has Nagata dimension $d$ with $d$-dimensional Nagata constant $\gamma$ then\footnote{We write $A\lesssim B$ to mean that $A\leq CB$ for some universal constant $C$. If $C$ depends on some parameter $\alpha$, we write $A\lesssim_\alpha B$ instead.}
\begin{equation}\label{naoreqintro}
    \SEP(M)\lesssim\gamma\cdot d^3.
\end{equation}
However, this bound is not completely satisfactory since the inequality that one recovers when considering a normed space in \eqref{naoreqintro} is considerably worse than the linear bound from \cite{Naor17}. Instead, we propose using the Lipschitz dimension introduced by Cheeger and Kleiner in \cite{CK13} as a more appropriate (or, at least, alternative) way to generalize \cite[Theorem 1.3]{Naor17} to arbitrary metric spaces. We obtain the following result.

\begin{theorem}\label{Lipdimsep}
    If $M$ is a metric space with finite Lipschitz dimension $d$ and $d$-dimensional Lipschitz distortion $\lambda$, then
     $$\SEP(M)\lesssim \lambda\cdot\sqrt{d}.$$
\end{theorem}

\noindent See Section \ref{sec:lipdim} for the definitions. Unlike \eqref{naoreqintro}, this upper bound agrees with the linear bound from \cite{Naor17} when $M$ is a normed space as $\lambda\lesssim\sqrt{d}$ in that case.

\subsection*{2. Finite Nagata dimension vs finite separation modulus}
Random decomposition moduli have been used in a wide variety of contexts and for a broad range of algorithmic and geometric applications: to bound distortion of stochastic embeddings \cite{Bar96}, to extend Lipschitz mappings \cite{LN05,Naor21} and to obtain embeddings into Hilbert spaces \cite{NS11} among others. Other works also considered such decompositions implicitly \cite{Lei88, AP90, LS91, KPR93, AKPW95, Rao99}. However, all of these examples involve metric spaces with finite Nagata dimension. More precisely, Bartal considers subsets of $\RR^d$ in \cite{Bar96}, Lee and Naor focus on doubling metric spaces in \cite{LN05,Naor21}, and Naor and Silberman stick to finite Nagata-dimensional spaces in \cite{NS11}. This hints towards the possibility that there might be a lower bound of some decomposition modulus in terms of the Nagata dimension, that is, a converse to inequality \eqref{naoreqintro}. Our next result shows that this is not the case.

\begin{theorem}\label{mainthm:nagatasep}
There exists a metric space $M$ with infinite Nagata dimension and finite separation modulus, which admits a stochastic embedding into trees with finite distortion.
\end{theorem}

\noindent Thus, finite Nagata dimension and finite separation modulus really are different properties.

\subsection*{3. Embedding Lipschitz-free spaces into $L_1$}

Embedding Lipschitz-free spaces (or equivalently, spaces of probability measures equipped with the Earth Mover's Distance) into $L_1$ is an important problem in image processing \cite{RTG00,NS07,AIK08,GO26} and the theory of normed spaces \cite{KN06,DKO20,BGS23}. In many results from the literature, the underlying metric space is required to have finite Nagata dimension (see, for instance, \cite{FG22,Gar25}). The natural question is then to what extent this condition is necessary, and the following result shows that it is artificial in some cases.

\begin{theorem}\label{mainthm:lipfree}
There exists a metric space $M$ with infinite Nagata dimension such that $\lipfree{M}$ is isomorphic to $\ell_1$. Moreover, there exists an increasing sequence $(S_n)$ of finite subsets of $M$ such that $M=\bigcup_nS_n$ and each $\lipfree{S_n}$ is $C$-isomorphic to $\ell_1^{\abs{S_n}-1}$ for a fixed $C<\infty$.
\end{theorem}

\noindent The second part of the statement provides an advance in the direction of the core problem raised in \cite[Problem 2.6]{DKO20} and studied in \cite{DKO20,DKO21,MT26}: namely, to obtain a metric parameter of finite metric spaces that measures the distance between their Lipschitz-free space and the corresponding $\ell_1^N$.

The example constructed in Theorem \ref{mainthm:lipfree} is stochastically embeddable into trees with finite distortion, although we will not prove it here. It is an open question, raised by Baudier and Schlumprecht (cf.~\cite[Question 1.2]{MT26}), whether there is a metric space that cannot be stochastically embedded into trees and whose Lipschitz-free space is isomorphic to $\ell_1$.

\subsection*{4. Separation modulus and stochastic embeddings into trees}

Theorem \ref{mainthm:lipfree} leads us to conclude that Nagata dimension is not the right dimensionality notion in the context of embeddings of Lipschitz-free spaces into $L_1$.
Alternatively, there are numerous works relating those embeddings to stochastic embeddings of the underlying metric space into trees (see \cite{MV23,BGS23,Gar25}; see also \cite{Bar96,CCGGP98,FRT04} where the alternative term ``random approximation by trees'' is used). The distortion of these stochastic embeddings of $M$, denoted here by $\STD(M)$ (see Definition \ref{def:std}), is indeed known to be closely related to the separation modulus since Bartal introduced the latter in \cite{Bar96} for this very reason. We develop the relation between $\SEP(M)$ and $\STD(M)$, complementing \cite[Theorem 13]{Bar96}, and show that both notions are equivalent for metric spaces with finite diameter and width.

\begin{theorem}\label{sepstdth}
   For every metric space $M$,
    $$\SEP(M)\lesssim \STD(M)\lesssim\SEP(M)\log\frac{2\diam(M)}{\omega(M)}$$
where $\omega(M)=\inf_{x\neq y\in M}d(x,y)$ and the right-hand side is understood to be $\infty$ if $\diam(M)=\infty$ or $\omega(M)=0$.
\end{theorem}

Both inequalities are tight. On the one hand, if $M$ is the $d$-dimensional grid of length $25$, then it is known that $\SEP(M)$ is of order $\Theta(d)$ \cite[Lemma 3.1]{CCGGP98}. Also, by \cite[Corollary 2]{KN06} and \cite[Corollary 19]{BMSZ22} it is known that $\STD(M)$ is of order $\Theta(d)$, matching the first inequality. On the other hand, if $M$ is the 2-dimensional grid of length $N>1$ then both
$\STD(M)$ and $\SEP(M)\log\frac{2\diam(M)}{\omega(M)}$ are of order $\Theta(\log N)$ (see \cite[Theorem 6.6]{AKPW95} and \cite{Naor17}).

The first inequality in Theorem \ref{sepstdth} is new, while the second is just a restatement of \cite[Theorem 13]{Bar96}. However, a variation of Bartal's argument allows us to obtain a stronger result for snowflaked metric spaces. Indeed, Gartland proved in \cite[Theorem 4.16]{Gar25} that snowflakes of finite Nagata-dimensional metric spaces are stochastically embeddable into trees. Our last result is a quantitative version of Gartland's theorem using the separation modulus.

\begin{theorem}\label{STDsnow}
    For every $\alpha\in(0,1)$ and every metric space $M$ we have
    $$\STD(M^\alpha)\lesssim_\alpha \SEP(M).$$
\end{theorem}

After Section \ref{sec:prelim} collecting common definitions, the remainder of this paper contains the proofs of Theorems \ref{Lipdimsep} to \ref{STDsnow}. The order in which they are proved will differ slightly from the order in which they appear in this introduction: we will prove Theorems \ref{sepstdth} and \ref{STDsnow} first as some of the other results depend on them.

\section{Preliminaries}
\label{sec:prelim}

In what follows, $(M,d)$ will always denote a metric space with metric $d$. The \emph{diameter} and \emph{width} of $M$ are defined as
$$
\diam(M)=\sup_{x,y\in M}d(x,y) \quad,\quad \omega(M)=\inf_{x\neq y\in M}d(x,y) ,
$$
respectively. Given $\alpha\in(0,1)$, $M^\alpha$ denotes the \textit{$\alpha$-snowflake} of $M$, that is, the metric space $(M,d^\alpha)$.

Some of the metric spaces considered here will be constructed as (undirected) graphs. The sets of vertices and edges of a graph $G$ will be denoted by $V(G)$ and $E(G)$, respectively. A \textit{weighted graph} is a graph where a positive weight is assigned to each edge. The \textit{geodesic distance} associated to such a graph is the length of the shortest path considering weights. We call \textit{graph distance} the geodesic distance obtained when assigning a weight of 1 to every edge. We are especially interested in (weighted) trees, that is, graphs without cycles which, when endowed with their geodesic distance, will be referred to as \textit{geodesic trees}.

We will sometimes assume that metric spaces are \textit{pointed}, that is, we have selected an arbitrary point $0\in M$ to act as base point.
We shall need the following construction, used to fuse different metric spaces together. Given a family $(M_\gamma,d_\gamma,0_\gamma)$ of pointed metric spaces, we define their \textit{metric sum} (cf. \cite[Definition 1.13]{Wea18}) as the disjoint union $M=\set{0}\cup\bigcup_\gamma(M_\gamma\setminus\set{0_\gamma})$ endowed with the metric $d$ given by
\begin{align*}
d(x,0) &= d_\gamma(x,0_\gamma) &&\text{if $x\in M_\gamma$} \\
d(x,x') &= d_\gamma(x,x') &&\text{if $x,x'\in M_\gamma$} \\
d(x,x') &= d_\gamma(x,0_\gamma) + d_{\gamma'}(x',0_{\gamma'}) &&\text{if $x\in M_\gamma,x'\in M_{\gamma'},\gamma\neq\gamma'$}
\end{align*}
Informally, $(M,d)$ is obtained by taking the disjoint union of the $M_\gamma$'s, collapsing all of their base points together, and having all paths between different $M_\gamma$'s go through the base point.

\medskip

Let $\mathcal{C}$ be a collection of subsets of $M$. Given $D>0$, we say that $\mathcal{C}$ is \textit{$D$-bounded} if each element of $\mathcal{C}$ has diameter at most $D$. Given $s>0$, the \textit{$s$-multiplicity} of $\mathcal{C}$ is the least number $n$ (possibly infinite) such that every subset of $M$ with diameter less than $s$ intersects at most $n$ elements of $\mathcal{C}$.

\begin{definition}[Nagata dimension]\label{def:nagata}
Given $d\in\NN\cup\set{0}$, the \textit{$d$-dimensional Nagata constant} of $M$, denoted $\gamma_d(M)$, is the smallest $\gamma>0$ (possibly infinite) with the property that for each $s>0$ there exists a $(\gamma s)$-bounded cover of $M$ with $s$-multiplicity at most $d+1$.

The \textit{Nagata dimension} of $M$, denoted $\dim_N(M)$, is the smallest $d$ such that $\gamma_d(M)<\infty$, or $\infty$ if there is no such $d$.
\end{definition}

It is clear that $\gamma_d(S)\leq\gamma_d(M)$ and $\dim_N(S)\leq\dim_N(M)$ for every $S\subset M$. Moreover, Nagata dimension is invariant under biLipschitz equivalence and under snowflaking \cite[Lemma 2.1]{LS05}. See \cite{LS05} for other basic facts on Nagata dimension.

\medskip

We will consider partitions $P$ of $M$; when the empty set is an element of $P$ (for instance, due to some of our constructions below), it will be disregarded. The unique cluster (i.e. element) of $P$ containing the point $x\in M$ will be denoted $P(x)$; we will also write $P(x)=\emptyset$ if $P$ is a partition of $S\subset M$ and $x\notin S$. We write $P\preccurlyeq P'$ if $P'$ is finer than $P$, that is, if $P'(x)\subset P(x)$ for every $x$. Given $S\subset M$, the restriction of $P$ to $S$, denoted $P\restrict_S$, is the partition of $S$ given by intersections of the clusters of $P$ with $S$, that is, $P\restrict_S=\set{A\cap S\,:\,A\in P}$.
Given $D>0$, the family of all $D$-bounded partitions of $M$ will be denoted $\mathcal{P}_D^M$. 

\begin{definition}[Separation modulus]\label{def:sep}
The \emph{separation modulus of $M$}, denoted $\SEP(M)$, is the infimum of $\lambda\in [1,\infty]$ such that, for every $D>0$ and every finite subset $S\subset M$, there exists a probability distribution $\mu$ on $\mathcal{P}_D^S$ such that
    \begin{equation}\label{eqsep}
        \mu\big(P\in\mathcal{P}_D^S\;:\;P(x)\neq P(y)\big)\leq \lambda\cdot\frac{d(x,y)}{D}\quad\forall x,y\in S.
    \end{equation}
\end{definition}

One can intuitively think of partitions on $M$ directly; the main reason why we consider partitions of $S$ instead is to avoid measurability problems in the definition of probability distributions on $\mathcal{P}_D^M$. See \cite[Section 1.7]{Naor21} for an equivalent definition using partitions of $M$.

\medskip

For the following definition, if $M$ is finite, we say that a metric $d_T$ on $M$ is a \emph{tree metric} if it agrees with the geodesic distance for some weighted tree $T$ with vertices in $M$. We say that $d_T$ \emph{dominates} $d$ if $d_T(x,y)\geq d(x,y)$ for all $x,y\in M$.

\begin{definition}[Stochastic tree distortion]\label{def:std}
The \emph{stochastic tree distortion of $M$}, denoted $\STD(M)$, is the infimum of $\rho\in [1,\infty]$ such that, for every finite subset $S\subset M$, there exists a probability distribution $p$ on the family of dominating tree metrics $d_T$ on $S$ for which
    \begin{equation}\label{eqstd}
        \EE_p(d_T(x,y)) \leq \rho\cdot d(x,y) \quad\forall x,y\in S.
    \end{equation}
\end{definition}

It is clear from the definitions that the separation modulus and the stochastic tree distortion are finitely determined, that is, $\SEP(M)=\sup\set{\SEP(S):S\subset M \text{ finite}}$ and similarly for $\STD$. Moreover, although they are defined as infima, if $M$ is finite then both are actually minima (i.e. the infima are attained) since we only need to consider probability distributions over a finite set. Indeed, this is clear for $\SEP(M)$ as there are only finitely many partitions of $M$. In the case of $\STD(M)$, there are only finitely many unweighted trees $T$ with set of vertices $M$, and for each one it suffices to consider the tree metric $d_T$ obtained by associating the weight $d(x,y)$ to each edge $\set{x,y}$, as any tree metric on $T$ dominating $d$ must dominate $d_T$ as well.

\section{Separation modulus and stochastic embeddings into trees}
\label{sec:stdsep}

We first focus on the proof of the first inequality of Theorem \ref{sepstdth}. For that purpose we introduce the following auxiliary results.

\begin{lemma}\label{subsep}
Let $M$ be finite and $\lambda\geq 1$. Suppose that for every $D>0$ there is a probability distribution $\mu$ on $\mathcal{P}_D^M$ such that
$$
\mu\big( P(x)\neq P(y) \big) \leq \lambda\cdot\frac{d(x,y)}{D} \qquad (x,y\in M).
$$
Then $\SEP(M)\leq\lambda$.
\end{lemma}

\begin{proof}
We only need to show that, for every $S\subset M$, there is a probability distribution $\mu_S$ on $\mathcal{P}_D^S$ satisfying \eqref{eqsep}. We may simply take
$$
\mu_S(P_S) = \mu\big(\set{P\in\mathcal{P}_D^M \,:\, P\restrict_S = P_S}\big) \qquad (P_S\in\mathcal{P}_D^S).
$$
    It is straighforward to check that $\mu_S$ is indeed a probability distribution (as every $P_S$ is the restriction of some $P$) and satisfies for every $x,y\in S$ that
    $$\mu_S\big(P_S(x)\neq P_S(y)\big)=\mu\big(P(x)\neq P(y)\big)\leq\frac{\lambda d(x,y)}{D}.$$
\end{proof}

\begin{lemma}\label{graphtree}
If $T$ is a finite unweighted tree and $d_T$ is its graph distance, then $\SEP(T,d_T)\leq 2$.
\end{lemma}

\begin{proof}
Fix a point $0\in T$ to use as its root. Given $x,y\in T$, let us write $x\preccurlyeq y$ if $x$ lies in the path between $0$ and $y$. For a given $x\in T$ and $r\geq 0$, we consider the set
$$
F(x,r) = \set{ y\in T \,:\, d_T(x,y)\leq r \text{ and } x\preccurlyeq y } .
$$
Note that this set has diameter at most $2r$.

Let $k\in\NN$. Given $m\in\set{1,\ldots,k+1}$, consider the family of sets
$$
P_m = \set{F(0,m-1)} \cup \set{ F(x,k) \,:\, x\in T \text{ and } d_T(x,0)\equiv m \text{ mod } k+1 } .
$$
It is readily checked that $P_m$ is a $2k$-bounded partition of $T$. Let $\mu_k$ be the probability distribution on $\mathcal{P}_{2k}^T$ given by $\mu_k(P_m)=\frac{1}{k+1}$ for each $m$, and $0$ for every other partition.
We claim that
$$
\mu_k\big(P(x)\neq P(y)\big) \leq \frac{d_T(x,y)}{k+1} \qquad x,y\in T .
$$
This is clear if $x=y$ or if the right-hand side is at least $1$, so we may assume $1\leq d_T(x,y)\leq k$. 
Let $z$ be the median point of $0$, $x$ and $y$, that is, the point in the path between $x$ and $y$ such that $z\preccurlyeq x$ and $z\preccurlyeq y$. Then $x,y\in F(z,d_T(x,y))$, and there are at least $k+1-d_T(x,y)$ values of $m$ for which $F(z,d_T(x,y))\subset P_m(z)$ and hence $P_m(x)=P_m(z)=P_m(y)$. Thus $P_m(x)\neq P_m(y)$ for at most $d_T(x,y)$ values of $m$, proving our claim.

The distributions $\mu_k$ witness that $\SEP(T)\leq 2$. Indeed, let $D>0$ and $x\neq y\in T$. If $D<2$, then for any $k\in\NN$ it is clear that
$$
\mu_k\big(P(x)\neq P(y)\big) \leq 1 < \frac{2}{D} \leq 2\,\frac{d_T(x,y)}{D} .
$$
Otherwise, let $k\in\NN$ be such that $2k\leq D<2(k+1)$. Then $\mu_k$ is a probability distribution on $\mathcal{P}_D^T\supset\mathcal{P}_{2k}^T$ and
$$
\mu_k\big(P(x)\neq P(y)\big) \leq \frac{d_T(x,y)}{k+1} < \frac{2(k+1)}{k+1} \frac{d_T(x,y)}{D} = 2\,\frac{d_T(x,y)}{D} .
$$
Thus $\SEP(T)\leq 2$ by Lemma \ref{subsep}.
\end{proof}

\begin{corollary}\label{treesep}
    If $(T,d_T)$ is a finite geodesic tree, then $\SEP(T)\leq 4$.
\end{corollary}
\begin{proof}

Construct a new geodesic tree $T'$ by subdividing the edges of $T$ so that each edge of $T'$ has weight between $\omega(T)$ and $2\omega(T)$. Then $T'$ contains an isometric copy of $T$ and is $2$-biLipschitz equivalent to a finite unweighted tree, thus $\SEP(T)\leq \SEP(T')\leq 4$ by Lemma \ref{graphtree}.		
\end{proof}

\begin{proof}[Proof of the first inequality of Theorem \ref{sepstdth}]

    We assume that $M$ is finite and take $D>0$, $\rho=\text{STD}(M)$. Consider a probability distribution $p$ over the family of dominating trees $\Omega$ witnessing stochastic distortion $\rho$. Then, for each $T\in \Omega$ we also consider by Corollary \ref{treesep} a probability distribution $\mu_T$ over $\mathcal{P}_D^T$ witnessing separation modulus $4$. Note that $\mathcal{P}_D^T\subset\mathcal{P}_D^M$ as $d_T\geq d$, so we may consider $\mu_T$ as a probability distribution over $\mathcal{P}_D^M$. Finally, we define a probability distribution $\mu$ over $\mathcal{P}_D^M$ by
    $$\mu(P)=\EE_p(\mu_T(P)) \qquad (P\in\mathcal{P}_D^M) .$$
    It is clear that $\mu$ is a probability distribution since
    $$\mu(\mathcal{P}_D^M) = \EE_p(\mu_T(\mathcal{P}_D^M)) = \EE_p(1) = 1.$$
    Also, if $x,y\in M$ are distinct then
		$$
		\mu\big(P(x)\neq P(y)\big) = \EE_p\big(\mu_T(P(x)\neq P(y))) \leq \EE_p\pare{\frac{4d_T(x,y)}{D}} \leq \frac{4\rho d(x,y)}{D} .
		$$
Thus $\mu$ witnesses $\SEP(M)\leq 4\rho$.
\end{proof}

Even though the second inequality of Theorem \ref{sepstdth} was already proven by Bartal in \cite[Theorem 13]{Bar96}, we include here a more detailed proof - not just for completeness, but also because an adjustment to Bartal's argument simultaneously yields the proof of Theorem \ref{STDsnow}. More specifically, we shall show for every $\alpha\in(0,1)$ that
\begin{equation}
    \STD(M^\alpha)\leq \frac{2^{\alpha+5}}{(2^\alpha-1)(2^{1-\alpha}-1)} \SEP(M).
\end{equation}

\begin{proof}[Proof of Theorem \ref{STDsnow} and the second inequality of Theorem \ref{sepstdth}]
    We may assume without loss of generality that $M$ is finite and $\omega(M)=1$. We denote $\Delta=\diam(M)$ and $\lambda=\SEP(M)$. Consider the least $t\in\NN$ such that $\frac{\Delta}{2^t}<1$. For each $i\in\{0,\dots,t\}$ and $S\subset M$ we denote $r_i:=\frac{\Delta}{2^i}$, $\mathcal{P}_i^S:=\mathcal{P}_{r_i}^S$ for $i\geq1$, $\mathcal{P}_0^S=\{\{S\}\}$ and
    $$\Omega_i^S=\{(P_j)_{j=i}^t\;:\;P_j\in\mathcal{P}_j^S\;\;,\;\;P_j\preccurlyeq P_{j+1}\}.$$
    Notice that $\mathcal{P}_t^S$ only contains the trivial partition whose clusters are the singletons of $S$. Also, for each $i\in\{0,\dots,t\}$ and $S\subset M$, we consider $\mu_i^S$ a probability distribution over $\mathcal{P}_i^S$ satisfying
    \begin{equation}\label{eqsep_both}
        \mu_i^S\big(P(x)\neq P(y)\big)\leq \lambda\cdot\frac{d(x,y)}{r_i}\quad\forall x,y\in S.
    \end{equation}
    
\begin{claim*}
There is a probability distribution $p$ over $\Omega_0^M$ such that
\begin{equation}\label{lamprop_bartalmain}
        p\big(P_k(x)\neq P_k(y)\big)\leq\frac{\lambda d(x,y)}{r_k}\;\;\;\;\;\;\forall x,y\in M\;\;\;\forall k\in\{0,\dots,t\}.
\end{equation}
\end{claim*}
In fact, we will prove that, for each $i\in\{0,\dots,t\}$ and $S\subset M$, there exists a probability distribution $p_i^S$ over the family $\Omega_i^S$ satisfying for every $x,y\in S$,
            \begin{equation}\label{lamprop_bartal}
        p_i^S\big(P_k(x)\neq P_k(y)\big)\leq\frac{\lambda d(x,y)}{r_k}\;\;\;\;\;\;(i\leq k\leq t).
    \end{equation}
Then $p$ will be chosen as $p_0^M$.

We prove it by induction on $i$. If $i=t$ then $\Omega_t^S$ consists of only one partition and $p_t^S$ is a Dirac delta which clearly satisfies \eqref{lamprop_bartal} since $\frac{\lambda d(x,y)}{r_t}>\lambda\geq1$ for $x\neq y$. Now, assuming that for some $i>0$ and every $S\subset M$, there exists a distribution $p_i^S$ satisfying \eqref{lamprop_bartal}, we define $p_{i-1}^S$ for each $(P_j)_{j=i-1}^t\in\Omega_{i-1}^S$ as
    $$p_{i-1}^S\big((P_j)_{j=i-1}^t\big)=\mu_{i-1}^S(P_{i-1})\prod_{R\in P_{i-1}}p_i^R\big((P_j\restrict_R)_{j=i}^t\big).$$
This probability may be understood as performing several random choices that determine a unique element in $\Omega_{i-1}^S$: first choose an initial partition $P_{i-1}\in\mathcal{P}_{i-1}^S$ with probability distribution $\mu_{i-1}^S$. Then, for each cluster $R\in P_{i-1}$, pick independently a sequence of increasing partitions of $R$ in $\Omega_i^R$ following the probabilistic distribution $p_i^R$, the union of which will give us an element of $\Omega_{i-1}^S$.
    
    We show now that $p_{i-1}^S$ is also a probability distribution satisfying \eqref{lamprop_bartal}. First notice that we have for all $S\subset M$ and all $P\in\mathcal{P}_{i-1}^S$
    \begin{equation}\label{star_bartal}
		\sum_{\substack{(P_j)_{j=i}^t\in\Omega_{i}^S\\P\preccurlyeq P_i}}\prod_{R\in P}p_i^R\big((P_j\restrict_R)_{j=i}^t\big) = \prod_{R\in P}\sum_{(P_j^R)_{j=i}^t\in\Omega_{i}^R}p_i^R\big((P^R_j)_{j=i}^t\big) = \prod_{R\in P}p_i^R(\Omega_i^R) = 1 .
		\end{equation}		
		So we get
    \begin{equation*}\begin{aligned}
		p_{i-1}^S(\Omega_{i-1}^S)=&\sum_{(P_j)_{j=i-1}^t\in\Omega_{i-1}^S}\mu_{i-1}^S(P_{i-1})\prod_{R\in P_{i-1}}p_i^R\big((P_j\restrict_R)_{j=i}^t\big)\\
		=&\sum_{P\in\mathcal{P}_{i-1}^S}\mu_{i-1}^S(P)\sum_{\substack{(P_j)_{j=i}^t\in\Omega_{i}^S\\P\preccurlyeq P_i}}\prod_{R\in P}p_i^R\big((P_j\restrict_R)_{j=i}^t\big)\\
		\stackrel{\eqref{star_bartal}}{=}&\sum_{P\in\mathcal{P}_{i-1}^S}\mu_{i-1}^S(P)=1.
		\end{aligned}\end{equation*}
It only remains to show that $p_{i-1}^S$ satisfies \eqref{lamprop_bartal}. We divide the proof into two different cases.
    
    \medskip
    \textbf{Case 1: $k=i-1$.} In this case we use the separation property \eqref{eqsep_both}:
    $$\begin{aligned}
		p_{i-1}^S(P_k(x)\neq P_k(y))=&\sum_{\substack{P_k\in\mathcal{P}_k^S\\P_k(x)\neq P_k(y) }}\mu_k^S(P_k)\sum_{\substack{(P_j)_j\in\Omega_{k+1}^S\\P_k\preccurlyeq P_{k+1} }}\prod_{R\in P_k}p_{k+1}^R\big((P_j\restrict_R)_j\big)\\
		\stackrel{\eqref{star_bartal}}{=}&\sum_{\substack{P_k\in\mathcal{P}_k^S\\P_k(x)\neq P_k(y) }}\mu_k^S(P_k)=\mu_k^S(P_k(x)\neq P_k(y))\stackrel{\eqref{eqsep_both}}{\leq}\frac{\lambda d(x,y)}{r_k}.
		\end{aligned}$$

    \medskip
    \textbf{Case 2: $k\geq i$.} In this case we use the induction hypothesis (IH) since
    $$\begin{aligned}
        p_{i-1}^S&(P_k(x)\neq P_k(y)) = \sum_{P_{i-1}\in\mathcal{P}_{i-1}^{S}}\mu_{i-1}^{S}(P_{i-1})\sum_{\substack{(P_j)_j\in\Omega_i^{S}\\P_{i-1}\preccurlyeq P_i\\P_k(x)\neq P_k(y) }}\prod_{R\in P_{i-1}}p_i^R\big((P_j\restrict_R)_j\big) \\
				=&\sum_{P_{i-1}\in\mathcal{P}_{i-1}^{S}}\mu_{i-1}^{S}(P_{i-1})\sum_{\substack{(P_j)_j\in\Omega_i^{P_{i-1}(x)}\\P_k(x)\neq P_k(y) }}p_{i}^{P_{i-1}(x)}\big((P_j)_j\big)\sum_{\substack{(P'_j)_j\in \Omega_i^{S\setminus P_{i-1}(x)}\\P_{i-1}\restrict_{S\setminus P_{i-1}(x)}\preccurlyeq P'_i }}\prod_{\substack{R\in P_{i-1}\\ R\neq P_{i-1}(x)}}p_i^R\big((P'_j\restrict_R)_j\big) \\
				\stackrel{\eqref{star_bartal}}{=}&\sum_{P_{i-1}\in\mathcal{P}_{i-1}^{S}}\mu_{i-1}^{S}(P_{i-1})\sum_{\substack{(P_j)_j\in\Omega_i^{P_{i-1}(x)}\\P_k(x)\neq P_k(y) }}p_{i}^{P_{i-1}(x)}\big((P_j)_j\big) \\
				=&\sum_{P_{i-1}\in\mathcal{P}_{i-1}^{S}}\mu_{i-1}^{S}(P_{i-1})p_i^{P_{i-1}(x)}\big(P_k(x)\neq P_k(y)\big) \stackrel{\mathrm{(IH)}}{\leq} \frac{\lambda d(x,y)}{r_k}
    \end{aligned}$$
	where we have used \eqref{star_bartal} with $S\setminus P_{i-1}(x)$ in place of $S$.
    This finishes the induction and thus our Claim is proven.

    \medskip
		
        Now fix $\alpha\in (0,1]$. For the remainder of the proof, we only consider the probability distribution $p$ and we rename its associated sample space $\Omega_0^M$ as $\Omega$. Every element $T=(P_j)_{j=0}^t\in\Omega$ defines a weighted tree whose vertices, edges and weights are given by
    $$\begin{aligned}&V(T)=\{(j,A)\;:\;j\in\{0,\dots,t\}\;,\;A\in P_j\},\\
    &E(T)=\big\{\{(j,A),(j-1,B)\}\;:\;j\in\{1,\dots,t\}\;,\;A\subset B\big\},\\
		&w\big(\{(j,A),(j-1,B)\}\big)=2^{\alpha-1}r_j^\alpha \quad \forall\{(j,A),(j-1,B)\}\in E(T).\end{aligned}$$	
We consider the geodesic distance $d_T$ in $V(T)$. By identifying each $x\in M$ with the leaf $(t,\{x\})\in V(T)$, $d_T$ becomes a metric on $M$. Now fix $x\neq y\in M$ and let $h_T$ be the largest integer $i\geq 0$ such that $P_i(x)=P_i(y)$; note that $h_T<t$. Then we have
\begin{equation}\label{treedist_both}
d_T(x,y)=2^\alpha\sum_{j=h_T+1}^tr_j^\alpha \qquad\text{and therefore}\qquad r_{h_T}^\alpha\leq d_T(x,y)\leq \frac{2^\alpha}{2^\alpha-1}r_{h_T}^\alpha
\end{equation}
as the path between $x$ and $y$ is the union of the paths between $x$ or $y$ and $P_{h_T}(x)=P_{h_T}(y)$.
In particular
$$
d(x,y)^\alpha \leq \diam(P_{h_T}(x))^\alpha \leq r_{h_T}^\alpha \leq d_T(x,y) ,
$$
so $d_T$ dominates $d^\alpha$. By a result of Gupta \cite[Theorem 1.1]{Gup18}, there exists a tree metric $d'_T$ on $M$ with $d_T\leq d'_T\leq 8d_T$. Then $p$ provides a probability distribution on the dominating tree metrics $d'_T$. Let us see that $p$ witnesses the desired bounds on $\STD(M^\alpha)$. We consider the cases $\alpha=1$ (Theorem \ref{sepstdth}) and $\alpha<1$ (Theorem \ref{STDsnow}) separately.

First, take $\alpha=1$. Then
    $$\begin{aligned}
        \mathbb{E}_p(d'_T(x,y))&=\sum_{T\in\Omega}p(T) d'_T(x,y)\stackrel{\eqref{treedist_both}}{\leq}\sum_{i=0}^{t-1}p(h_T=i)\cdot 16r_i \\
				&=\sum_{i=0}^{t-1}p\big(P_i(x)=P_i(y)\;\text{ and }\;P_{i+1}(x)\neq P_{i+1}(y)\big)\cdot 16r_i \\
				&\leq\sum_{i=0}^{t-1}p\big(P_{i+1}(x)\neq P_{i+1}(y)\big)\cdot 16r_i \\
				&\stackrel{\eqref{lamprop_bartalmain}}{\leq}\sum_{i=0}^{t-1}\frac{\lambda d(x,y)}{r_{i+1}}\cdot 16r_i\leq 32t\lambda d(x,y).
    \end{aligned}$$
Thus $\STD(M)\leq 32t\lambda$ with $t\leq 1+\log_2\Delta$. This proves the second inequality of Theorem \ref{sepstdth}.

Now suppose $\alpha<1$. Let $h\in\{0,1,\ldots,t-1\}$ be such that $r_{h+1}<d(x,y)\leq r_h$. Then
    $$\begin{aligned}
        \mathbb{E}_p(d'_T(x,y))&=\sum_{T\in\Omega}p(T) d'_T(x,y)\stackrel{\eqref{treedist_both}}{\leq}\sum_{i=0}^{t-1}p(h_T=i)\cdot 8\cdot 2^\alpha\sum_{j=i+1}^tr_j^\alpha \\
				&=2^{\alpha+3} \sum_{i=0}^{t-1}p\big(P_i(x)=P_i(y)\;\text{ and }\;P_{i+1}(x)\neq P_{i+1}(y)\big) \sum_{j=i+1}^tr_j^\alpha \\
				&\leq 2^{\alpha+3} \sum_{i=0}^h p\big(P_{i+1}(x)\neq P_{i+1}(y)\big) \sum_{j=i+1}^tr_j^\alpha \\
				&\stackrel{\eqref{lamprop_bartalmain}}{\leq} 2^{\alpha+3} \sum_{i=0}^h\frac{\lambda d(x,y)}{r_{i+1}} \sum_{j=i+1}^tr_j^\alpha \\
				&\leq 2^{\alpha+3}\lambda \sum_{i=0}^h\frac{r_h}{r_{i+1}} \sum_{j=i+1}^tr_j^\alpha
				\leq 2^{\alpha+3}\lambda \sum_{i=0}^h 2^{i+1-h} \frac{r_i^\alpha}{2^\alpha-1} .
    \end{aligned}$$
Since
$$
\sum_{i=0}^h 2^ir_i^\alpha = \Delta^\alpha \sum_{i=0}^h 2^{-i(1-\alpha)} \leq \Delta^\alpha \frac{2^{(1-\alpha)(h+1)}}{2^{1-\alpha}-1} = r_{h+1}^\alpha \frac{2^{h+1}}{2^{1-\alpha}-1} ,
$$
we conclude
$$
\mathbb{E}_p(d'_T(x,y)) \leq 2^{\alpha+3}\lambda \frac{2^{1-h}}{2^\alpha-1} r_{h+1}^\alpha \frac{2^{h+1}}{2^{1-\alpha}-1} \leq \frac{2^{\alpha+5}}{(2^\alpha-1)(2^{1-\alpha}-1)} \lambda d(x,y)^\alpha
$$
which proves Theorem \ref{STDsnow}.
\end{proof}

\section{Separation modulus and Lipschitz dimension}
\label{sec:lipdim}

In this section, we will prove Theorem \ref{Lipdimsep}.
We start by introducing the necessary notions, which are due to Cheeger and Kleiner in \cite{CK13}.

Given $r>0$, a finite sequence $(x_1,\dots,x_k)$ in $M$ is an \textit{$r$-path} if $d(x_i,x_{i+1})\leq r$ for every $i\in\{1,\dots,k-1\}$.
Two points $x,y$ in $M$ are \textit{$r$-connected} if there is an $r$-path $(x_1,\dots,x_k)$ between them, that is, $x_1=x$ and $x_k=y$. This defines an equivalence relation whose equivalence classes are called the \textit{$r$-components} of $M$.

    A map $f:M\to N$ between metric spaces is \textit{Lipschitz light} if there is a constant $C>0$ such that $f$ is Lipschitz with constant $C$ and, for every $r>0$ and every $W\subset N$ with $\diam(W)\leq r$, the $r$-components of $f^{-1}(W)$ have diameter at most $Cr$.   
    Equivalently, there are constants $c_1,c_2\geq1$ such that every $r$-connected set $K\subset M$ satisfies
    \begin{equation}\label{diamLipdim}c_1^{-1}\,\diam(K)-r\leq \diam(f(K))\leq c_2\,\diam(K).\end{equation}
    In that case, we say that $f$ is \textit{$(c_1,c_2)$-Lipschitz light}. The least possible value of $c_1 c_2$ satisfying the latter inequalities is known as the \textit{Lipschitz light distortion} of $f$, denoted here as $c(f)$.

\begin{definition}[Lipschitz dimension]\label{def:lipdim}
    The \textit{Lipschitz dimension} of a metric space $M$, denoted $\dim_L(M)$, is the least $d\in\NN$ such that there is a Lipschitz light map $f:M\to\ell_2^d$. If $M$ admits no Lipschitz light map into any Euclidean space, we put $\dim_L(M)=\infty$.

    If $\dim_L(M)=d<\infty$, we define the \textit{$d$-dimensional Lipschitz distortion} of $M$ as 
    \begin{equation}\label{eq:ldd}
    \lambda=\inf\{c(f)\;:\; f:M\to\ell_2^d\;\text{ Lipschitz light}\}.
    \end{equation}
\end{definition}

We have $\dim_N(M)\leq\dim_L(M)$ for any metric space $M$ by \cite[Corollary 3.5]{Dav19}. On the other hand, there exist metric spaces with finite Nagata dimension and infinite Lipschitz dimension, for instance non-abelian Carnot groups \cite[Theorem 6.8]{Dav19}.

\begin{proof}[Proof of Theorem \ref{Lipdimsep}]
    Set $d=\text{dim}_L(M)$, and let $f:M\to\ell_2^d$ be a $(c_1,c_2)$-Lipschitz light map with $c_1,c_2\geq 1$.
    Fix $D>0$ and a finite subset $S\subset M$, and set $r=\frac{D}{2c_1}$.
    For every partition $P\in\mathcal{P}_r^{f(S)}$, we define the partition $P_f\in\mathcal{P}_{D}^S$ whose clusters are the $r$-components of the inverse image of clusters of $P$ by $f\restrict_S$, that is,
    $$P_f=\bigg\{K\subset S\;:\;K\;\text{ is a $r$-component of }f^{-1}(A)\cap S\text{ for some }A\in P\bigg\}.$$
    It is straighforward to check that $P_f$ is really a partition of $S$ and that, for every $K\in P_f$, we have $\text{diam}(K)\leq D$ by \eqref{diamLipdim}. Now, since $\SEP(\ell_2^d)\leq 2\sqrt{d}$ by \cite[Theorem 3.1]{CCGGP98}, there exists a probability distribution $\mu$ over $\mathcal{P}_r^{f(S)}$ with $r$-separation constant $2\sqrt{d}$. We define now a probability distribution $\mu_f$ over $\mathcal{P}_{D}^S$ as
    $$\mu_f(Q)=\mu\big(\set{P\in\mathcal{P}_r^{f(S)} : P_f=Q}\big) \qquad(Q\in\mathcal{P}_{D}^S).$$
    It is clear that $\mu_f(\mathcal{P}_D^S)=\mu(\mathcal{P}_r^{f(S)})=1$.
    Let us show that $\mu_f$ has a $D$-separation constant bounded by $4c_1c_2\sqrt{d}$. Indeed, if we consider $x,y\in S$ with $d(x,y)\leq r$ then $x$ and $y$ belong to the same $r$-component of $S$ so that, for every $P\in\mathcal{P}_r^{f(S)}$,
    \begin{equation}\label{equivpart}
        P_f(x)\neq P_f(y)\;\;\;\Longrightarrow\;\;\;P(f(x))\neq P(f(y)).
    \end{equation}
    Therefore,
		\begin{align*}
		\mu_f\big(Q(x)\neq Q(y)\big) &= \mu\big(P_f(x)\neq P_f(y)\big) \\
		&\stackrel{\eqref{equivpart}}{\leq} \mu\big(P(f(x))\neq P(f(y))\big) \\
		&\leq 2\sqrt{d}\,\frac{d(f(x),f(y))}{r} \leq 4c_1c_2\sqrt{d}\,\frac{d(x,y)}{D} .
		\end{align*}
    On the other hand, if $d(x,y)\geq r$ then
    $$\mu_f\big(Q(x)\neq Q(y)\big)\leq 1\leq\frac{4c_1c_2\sqrt{d}\,d(x,y)}{D}.$$
    Thus $\SEP(M)\leq 4c_1c_2\sqrt{d}$, and taking the infimum over $f$ yields the desired conclusion.
\end{proof}

It is important to mention that the bound from Theorem \ref{Lipdimsep} is tight: if $M=\ell_1^d$ then $\text{dim}_L(M)=d$ with $\lambda\leq\sqrt{d}$ so that the right hand side of Theorem \ref{Lipdimsep} is upper bounded by $4d$, but the left hand side is known to be of order $\Theta(d)$ \cite[Theorem 3.3]{CCGGP98}. This argument yields the following conclusion.

\begin{corollary}
    For every $d\in\NN$, the $d$-dimensional Lipschitz distortion of $\ell_1^d$ is of order $\Theta(\sqrt d)$.
\end{corollary}

This means that one cannot get better distortion of Lipschitz light maps from $\ell_1^d$ into $\ell_2^d$ than the one obtained via linear isomorphisms.

In view of \cite[Theorem 1.5]{Dav23} and \cite[Corollary 3.5]{Dav19}, the family of snowflakes of spaces with finite Lipschitz dimension agrees with the family of snowflakes of spaces with finite Nagata dimension. Therefore, the following consequence of Theorem \ref{STDsnow} and Theorem \ref{Lipdimsep} may be understood as a quantitative version of \cite[Theorem 4.16]{Gar25}.

\begin{corollary}\label{quantgart}
    If $M$ is a metric space with finite Lipschitz dimension $d$ and $d$-dimensional Lipschitz distortion $\lambda$ then, for every $\alpha\in (0,1)$,
     $$\STD(M^\alpha)\lesssim_\alpha \lambda\cdot\sqrt{d}.$$
\end{corollary}

\section{Separation modulus and Nagata dimension: An example}\label{ex1}

In this section, we will construct metric spaces with arbitrarily high Nagata dimension whose stochastic tree distortion is uniformly bounded, and combine them to obtain a metric space $M$ with $\dim_N(M)=\infty$ and $\STD(M)<\infty$ which, by Theorem \ref{sepstdth}, proves Theorem \ref{mainthm:nagatasep}.

To that end, we consider the grid $\ZZ^d$ for $d\geq 2$, endowed with two different metrics. On one hand, the $\infty$-metric induced by the supremum norm $\norm{\cdot}_\infty$ on $\RR^d$. On the other hand, the distance $\rho$ given by
$$
\rho(x,y) = \begin{cases}
1+\log_{d^2}\norm{x-y}_\infty &\text{, if $x\neq y$} \\
0 &\text{, if $x=y$}
\end{cases} .
$$
It is straightforward to check that $\rho$ is a metric for $d\geq 2$. We denote by $B_\infty(x,r)$ and $B_\rho(x,r)$ the closed balls with center $x$ and radius $r$ with respect to the $\infty$-metric and $\rho$, respectively. Note that we have
    \begin{equation}
        B_\rho(x,r)=B_\infty(x,d^{2(r-1)}).
    \end{equation}

We will need the following well-known combinatorial lemma. For every $N,d\in\NN$ we denote $[N]^d=\{0,\dots,N\}^d$ and, for $i\in\{1,\dots,d\}$ and $x=(x_1,\dots,x_d)\in[N]^d$ we denote $x(i)=x_i$.

\begin{lemma}[Cubic version of Sperner's Lemma, \cite{Kuh60}]\label{sperner}
    Given $N,d\in\NN$, suppose that a coloring $f:[N]^d\to\{1,\dots,d+1\}$ satisfies the following two conditions:
    \begin{itemize}
        \item If $x(i)=0$ then $f(x)\leq i$.
        \item If $x(i)=N$ then $f(x)\neq i$.
    \end{itemize}
    Then, there is a unit cube $C$ of $[N]^d$ such that $f(C)=\{1,\dots,d+1\}$.
\end{lemma}

\begin{proposition}\label{pr14}
For each $N\in\NN$, the $(d-1)$-dimensional Nagata constant of $([N]^d,\rho)$ is at least $\log_{d^2}(N)$. Consequently, $\dim_N(\ZZ^d,\rho)\geq d$.
\end{proposition}

\begin{proof}
Let $\mathcal{C}$ be a $\log_{d^2}(N)$-covering of $([N]^d,\rho)$. It will be enough to show that the $1$-multiplicity of $\mathcal{C}$ is at least $d+1$.
For $i=1,\ldots,d$ let $A_i$ be the set of points $x\in [N]^d$ whose $i$-th coordinate is $0$. For every $x\in[N]^d$ denote by $\mathcal{C}(x)$ the family of clusters of $\mathcal{C}$ containing $x$. Consider then a coloring
    $f:[N]^d\to\{1,\dots,d+1\}$ that assigns to each $x\in[N]^d$
    $$f(x)=\begin{cases}
        \min\big\{i=1,\dots,d\;:\;\exists C\in\mathcal{C}(x),\;C\cap A_i\neq\emptyset\big\},\\
        d+1\quad\text{if the set above is empty.}
    \end{cases}$$
		It is straightforward to check that $f$ satisfies the hypotheses of Lemma \ref{sperner}. Therefore, there exists a unit cube $C$ in $[N]^d$ satisfying that $f(C)=\{1,\dots,d+1\}$. One easily checks that vertices $x\in C$ having different values of $f(x)$ must belong to different clusters of $\mathcal{C}$, and thus the multiplicity of $C$ is at least $d+1$. Our claim then follows because $\diam_\rho(C)=1$.
\end{proof}

We will now show that $\STD(\ZZ^d,\rho)$ is uniformly bounded. In fact, we will prove something a bit more general. Let us recall another modulus essentially introduced in \cite{NS11}, alternative to Definition \ref{def:sep}, that deals with a different type of random decomposition.

\begin{definition}[Padded decomposition modulus]\label{def:pad}
Given $\delta\in(0,1)$, the \emph{$\delta$-padded decomposition modulus of $M$} is the infimum of $\lambda\in [1,\infty]$ such that, for every $D>0$ and every finite subset $S\subset M$, there exists a probability distribution $\mu$ on $\mathcal{P}_D^S$ such that
    \begin{equation}\label{eqpad}
        \mu\pare{ P\in\mathcal{P}_D^S\;:\;B\pare{x,\frac{D}{\lambda}}\subset P(x) } \geq\delta \quad\forall x\in S.
    \end{equation}
\end{definition}

This definition can be found e.g. in \cite{Naor17} for finite metric spaces and \cite{LN25} (or, more rigorously, in \cite[Definition 3.2]{JPPQ26}) for infinite metric spaces. Again, we use the finite formulation to avoid measurability considerations, which are taken into account in \cite{JPPQ26}. By the proof of \cite[Theorem 2.2]{LN25} (see \cite[after Definition 2]{Naor17}), if the $\delta$-padded decomposition modulus of $M$ is $\lambda$ then $\SEP(M)\leq 4\lambda/\delta$.

Analogously to the separation modulus (Lemma \ref{subsep}), if $M$ is finite then, in order to compute the padded decomposition modulus, it is enough to consider probability distributions satisfying \eqref{eqpad} only on partitions of the whole space $M$ instead of its proper subsets.

\begin{proposition}\label{Hyppad}
    The $\frac{1}{2}$-padded decomposition modulus of $(\ZZ^d,\rho)$ is at most $9$.
\end{proposition}
\begin{proof}
    Fix $D>0$ and $N\in\NN$, and denote the subset $[-N,N]^d=\{-N,\dots,N\}^d\subset\ZZ^d$ and $M_N=([-N,N]^d,\rho)$.
		We will construct a probability distribution $p$ on $\mathcal{P}_D^{M_N}$ satisfying
    \begin{equation}\label{eqclaimpro}
        p\pare{M_N\cap B_{\rho}\pare{x,\frac{D}{\lambda}}\subset P(x)}\geq\frac{1}{2}
    \end{equation}
		with $\lambda=9$. This will be enough by the comment above. If $D<8$ then $B_{\rho}(x,D/9)=\set{x}$, so \eqref{eqclaimpro} is trivial. We may therefore assume $D\geq 8$.
		
		Suppose first that $D$ is an integer. We will prove that \eqref{eqclaimpro} holds with $\lambda=8$. We define $p$ as follows. Let $M'_N=[-(N+d^D),N+d^D]^d$ and consider the family $\Sigma_N$ of enumerations of $M'_N$, that is, injective functions $\sigma:\set{1,\ldots,\abs{M'_N}}\to M'_N$. For every $\sigma\in\Sigma_N$ we generate a partition $P_\sigma$ of $M_N$ as follows:
    $$P_\sigma=\set{ M_N\cap B_{\rho}\Big(\sigma(k),\frac{D}{2}\Big) \setminus\bigcup_{j=1}^{k-1} B_{\rho}\Big(\sigma(j),\frac{D}{2}\Big) }_{k=1}^{|M'_N|}.$$
    It is clear that $P_\sigma$ belongs to $\mathcal{P}_D^{M_N}$ for every $\sigma\in\Sigma_N$. Now, let $\mu$ be the uniform probability over $\Sigma_N$. We define a probability distribution $p_D$ on $\mathcal{P}_D^{M_N}$ as
    $$p_D(P)=\mu\pare{\set{\sigma\in\Sigma_{N} \,:\, P_\sigma=P}} \qquad(P\in\mathcal{P}_D^{M_N}).$$
    
		Let us show that $p_D$ satisfies \eqref{eqclaimpro} with $\lambda=8$. For that purpose, we consider, for every $\sigma\in \Sigma_{N}$ and $x\in M_N$, the index
    $$k_{\sigma,x}=\min\{k\in\NN\;:\; \sigma(k)\in B_{\infty}(x,d^{D-2}+d^{D/4-2})\}.$$
		We claim the following
		\begin{equation}\label{claim_Hyppad}
		\text{If $\sigma(k_{\sigma,x})\in B_\infty(x,d^{D-2}-d^{D/4-2})$, then $M_N\cap B_\rho(x,D/8)\subset P_\sigma(x)$.}
		\end{equation}
		Indeed, if that holds and $y\in M_N$ with $y\in B_\rho(x,D/8)=B_{\infty}(x,d^{D/4-2})$, then
    $$k_{\sigma,x}=\min\big\{k\in\NN\;:\; \sigma(k)\in B_\rho(y,D/2)\big\},$$
    because $B_\infty(x,d^{D-2}+d^{D/4-2})\supset B_\rho(y,D/2)=B_{\infty}(y,d^{D-2})\supset B_\infty(x,d^{D-2}-d^{D/4-2})$. Therefore,
    $$P_\sigma(y)=M_N\cap B_{\rho}\Big(\sigma(k_{\sigma,x}),\frac{D}{2}\Big) \setminus\bigcup_{j=1}^{k_{\sigma,x}-1} B_{\rho}\Big(\sigma(j),\frac{D}{2}\Big) =P_\sigma(x).$$
This proves \eqref{claim_Hyppad}. Using our claim and the fact that $\mu$ is uniform, we have that
    $$\begin{aligned}
		p_D\big(M_N&\cap B_{\rho}(x,D/8)\subset P(x)\big) = \mu\big(M_N\cap B_{\rho}(x,D/8)\subset P_\sigma(x)\big) \\
		&\stackrel{\eqref{claim_Hyppad}}{\geq} \mu\big(\sigma(k_{\sigma,x})\in B_\infty(x,d^{D-2}-d^{D/4-2})\big)\\
		&=\frac{|B_\infty(x,d^{D-2}-d^{D/4-2})|}{|B_\infty(x,d^{D-2}+d^{D/4-2})|}=\frac{(2(d^{D-2}-d^{D/4-2})+1)^d}{(2(d^{D-2}+d^{D/4-2})+1)^d}\geq\bigg(\frac{d^{D}-d^{D/4}}{d^{D}+d^{D/4}}\bigg)^d\\
		&\geq\bigg(\frac{d^{D}-2d^{D/4}}{d^{D}}\bigg)^d\geq\bigg(\frac{d^{D/2}-1}{d^{D/2}}\bigg)^d.\end{aligned}$$
    It now suffices to show that $(d^{D/2}-1)^d\geq\frac{1}{2}d^{dD/2}$. This holds since
    \begin{align*}
        (d^{D/2}-1)^d &\geq d^{dD/2}-\sum_{k=0}^{d-1}\binom{d}{k}(d^{D/2})^k \geq d^{dD/2}-\sum_{k=0}^{d-1}d^{d-k}(d^{D/2})^k \\
				&= d^{dD/2}-d^d\sum_{k=0}^{d-1}(d^{D/2-1})^k \geq d^{dD/2}-2d^{dD/2+1-D/2} \geq \frac{1}{2}d^{dD/2}
    \end{align*}
		where we use $d\geq 2$ and $D\geq 8$.
		
		Finally, if $D>8$ is not an integer, then let $K$ be the largest integer smaller than $D$. Then we may consider $p_K$ as a probability distribution on $\mathcal{P}_D^{M_N}\supset\mathcal{P}_K^{M_N}$, and we have
		$$
		p_K\pare{M_N\cap B_{\rho}\pare{x,\frac{D}{9}}\subset P(x)} \geq p_K\pare{M_N\cap B_{\rho}\pare{x,\frac{K}{8}}\subset P(x)} \geq \frac{1}{2}
		$$
		as $D/9\leq K/8$. This ends the proof.
\end{proof}

As a consequence of Proposition \ref{Hyppad} and \cite[Theorem 2.2]{LN25}, we get the following property of $(\ZZ^d,\rho)$.

\begin{corollary}\label{sep64}
    The separation modulus of $(\ZZ^d,\rho)$ is uniformly bounded for all $d\geq 2$.
\end{corollary}

We may now apply our results from Section \ref{sec:stdsep} to obtain spaces with uniformly bounded stochastic tree distortion.

\begin{corollary}\label{cor19}
    For every $d\geq 2$, there is a metric space $M_d$ with $\dim_N(M_d)\geq d$ and $\STD(M_d)\lesssim1$.
\end{corollary}
\begin{proof}
    Fix $\alpha\in (0,1)$ and let $M_d$ be the $\alpha$-snowflake of $(\ZZ^d,\rho)$. Then $\dim_N(M_d)=\dim_N(\ZZ^d,\rho)$ by \cite[Lemma 2.1]{LS05}. The statement thus follows from Theorem \ref{STDsnow} and Corollary \ref{sep64}.
\end{proof}

Finally, we will combine the spaces $M_d$ into a single one. We need a couple of auxiliary lemmas. The first one is a counterpart of Lemma \ref{subsep} for stochastic tree distortion.

\begin{lemma}\label{substd}
Let $M$ be finite and $\lambda\geq 1$. Suppose that there is a probability distribution $p$ on the family of dominating tree metrics $\tau$ on $M$ such that
$$
\EE_p(\tau(x,y)) \leq \lambda\cdot d(x,y) \quad\forall x,y\in M.
$$
Then $\STD(M)\leq 8\lambda$.
\end{lemma}

\begin{proof}
Let $S\subset M$. Given any dominating tree metric $\tau$ on $M$, there is by Gupta's theorem \cite[Theorem 1.1]{Gup18} a dominating tree metric $d_\tau$ on $S$ such that $d_\tau(x,y)\leq 8\tau(x,y)$ for $x,y\in S$. Thus we can define a probability distribution $\mu$ on the dominating tree metrics on $S$ by
$$
\mu(\rho)=p(\set{\tau\,:\,\rho=d_\tau})
$$
and it is clear that $\EE_\mu(\rho(x,y))=\EE_p(d_\tau(x,y))\leq 8\EE_p(\tau(x,y))\leq 8\lambda d(x,y)$ for $x,y\in S$.
\end{proof}

\begin{lemma}\label{sumstd}
Let $(M_\gamma,d_\gamma,0_\gamma)$ be a family of pointed metric spaces and let $(M,d,0)$ be their metric sum. Then $\STD(M)\leq 8\sup_\gamma\STD(M_\gamma)$.
\end{lemma}

\begin{proof}
Set $C=\sup_\gamma\STD(M_\gamma)$, and let $S$ be a non-empty finite subset of $M$ containing $0$. We will show that there is a probability distribution $p$ on the dominating tree metrics $\tau$ on $S$ such that $\EE_p(\tau(x,y))\leq Cd(x,y)$ for $x,y\in S$. The result will then follow from Lemma \ref{substd} as any finite subset of $M$ is contained in one such $S$, simply by adding the base point.

We may write $S=S_1\cup\ldots\cup S_n$ where, for $1\leq i\leq n$, $S_i$ is a finite non-empty subset of some $M_{\gamma_i}$ that contains $0$. Suppose that we are given dominating tree metrics $\tau_i$ on $S_i$ for $1\leq i\leq n$, and consider the metric $\tau = \tau(\tau_1,\ldots,\tau_n)$ on $S$ such that $(S,\tau)$ is the metric sum of the sets $(S_i,\tau_i,0)$; that is, $\tau(x,y)=\tau_i(x,y)$ for $x,y\in S_i$ and $\tau(x,y)=\tau_i(x,0)+\tau_j(0,y)$ for $x\in S_i$, $y\in S_j$. 
This is a tree metric on $S$, where the underlying tree structure contains precisely all edges between points of $S_i$ given by $\tau_i$ for each $i$, and it dominates $d$: this is clear for pairs of points in $S_i$ as $\tau$ equals $\tau_i$ there, and for $x\in S_i$, $y\in S_j$ we have
$$
\tau(x,y) = \tau_i(x,0)+\tau_j(0,y) \geq d(x,0)+d(0,y) = d(x,y) .
$$

By assumption, for each $i$ there exists a probability distribution $p_i$ on the dominating tree metrics $\tau_i$ on $S_i$ such that $\EE_{p_i}(\tau_i(x,y))\leq Cd(x,y)$. Then we define the probability distribution $p$ on the dominating tree metrics on $S$ by
$$
p(\tau) = \sum\set{ \prod_{i=1}^n p_i(\tau_i) \,:\, \tau=\tau(\tau_1,\ldots,\tau_n) } .
$$
For $x,y\in S_i$, we have $\EE_p(\tau(x,y)) = \EE_{p_i}(\tau_i(x,y)) \leq Cd(x,y)$, and for $x\in S_i$, $y\in S_j$,
$$
\EE_p(\tau(x,y)) = \EE_p(\tau_i(x,0)) + \EE_p(\tau_j(0,y)) \leq Cd(x,0) + Cd(0,y) = Cd(x,y) .
$$
This ends the proof.
\end{proof}

\begin{proof}[Proof of Theorem \ref{mainthm:nagatasep}]
Fix arbitrary base points $0_d\in M_d$ for every $d\geq 2$ and let $M$ be the metric sum of the spaces $M_d$. Then $\dim_N(M)=\infty$ because $M$ contains isometric copies of all $M_d$, and $\STD(M)<\infty$ by Corollary \ref{cor19} and Lemma \ref{sumstd}.
\end{proof}

\section{Nagata dimension and $\ell_1$-distortion: An example}

Finally, we will prove Theorem \ref{mainthm:lipfree} by constructing, for every $d\geq 6$, a metric space $(M_d,\rho_d)$ with Nagata dimension at least $d$ and satisfying, for some universal $C>0$ independent of $d$, that
\begin{equation}\label{dBMeq}d_{\mathrm{BM}}(\lipfree{M_d},\ell_1)\leq C,\end{equation}
where $d_\mathrm{BM}$ denotes Banach-Mazur distance. We will then combine all $M_d$ together to obtain a metric space $M$ with infinite Nagata dimension whose Lipschitz-free space is isomorphic to $\ell_1$.

While we will have no need to deal with Lipschitz-free spaces directly, we briefly recall their definition here. Given a pointed metric space $M$, let $\Lip_0(M)$ be the Banach space of all Lipschitz functions $f:M\to\RR$ such that $f(0)=0$, endowed with the norm given by the optimal Lipschitz constant. The \textit{Lipschitz-free space} $\lipfree{M}$ is then the canonical predual of $\Lip_0(M)$, which can be identified with the closed subspace of $\Lip_0(M)^*$ generated by all evaluation functionals $f\mapsto f(x)$ for $x\in M$. See \cite{GK03,Wea18} for reference on these spaces, and \cite{KN06,DKO20,MT26} for an alternative construction of $\lipfree{M}$ as the \textit{transportation cost space} over $M$.

We now proceed with our construction. For each $d$, the space $(M_d,\rho_d)$ will be a hyperbolic approximation of the flat torus $(\TT^d,\|\cdot\|_\infty)$ for a specific choice of hyperbolic parameters (see \cite{BS07, MT26}). For the sake of completeness, we include here a self-contained description of the space. We first construct a graph $G_d$ whose set of vertices, when endowed with the graph distance, will constitute our sought space $M_d$. For that purpose, we fix $r=1/d^2$ and, for $i\in\ZZ$, we let $L_i$ be the $r^i$-grid in $(\TT^d,\|\cdot\|_\infty)$, that is, $L_i=\big(r^i\ZZ^d\big)\cap\TT^d$. It is clear that $L_i=\{0\}$ for every $i<0$. Henceforth, $B_\infty(x,R)$ will denote the closed ball in $(\TT^d,\|\cdot\|_{\infty})$ of center $x\in\TT^d$ and radius $R>0$.

Firstly, we set the vertices of $G_d$ as $V(G_d)=\bigsqcup_{i=-\infty}^{\infty}L_i$ where, for every $i\in\ZZ$, we will refer to $L_i$ as the $i^\text{th}$ layer of $G_d$. There will be two types of edges in $G_d$.
\begin{itemize}
\item \textit{Horizontal edges} between elements from the same layer. Two vertices $x,y\in L_i$ are connected by an edge if $B_\infty(x,r^{i-1})\cap B_\infty(y, r^{i-1})\neq \emptyset$; that is, if $\norm{x-y}_\infty\leq 2r^{i-1}$.
\item \textit{Radial edges} between elements from consecutive layers. Vertices $x\in L_i$ and $y\in L_{i+1}$ are connected by an edge if $B_\infty(y, r^{i})\subset B_\infty(x,r^{i-1})$; that is, if $\norm{x-y}_\infty\leq r^{i-1}-r^i$.
\end{itemize}

We will prove \eqref{dBMeq} by finding an increasing sequence $(M^n_d)_{n\in\NN}$ of finite metric subspaces of $M_d$ such that $\bigcup_{n\in\NN}M_d^n=M_d$ and $d_{\mathrm{BM}}(\mathcal{F}(M_d^n),\ell_1^N)\leq C$ for universal $C>0$ and $N=|M_d^n|-1$. In fact, for every $n\in\NN$, the space $M_d^n$ will be the subset of $M_d$ whose points are $\bigsqcup_{i=-n}^{n}L_i$.

It is important to mention that $M_d^n$ is the metric space that arises when one considers the subgraph of $G_d$ given by the vertices $\bigsqcup_{i=-n}^{n}L_i$ with the graph distance (see \cite[Lemma 6.2.5]{BS07}).

Also, it is worth mentioning that this hyperbolic approximation matches those of \cite{BS07, MT26} when considering $\lambda=1/r$ there. More precisely, $M_d^n$ is the $n^\text{th}$ element of the sequence of $(\lambda,r)$-hyperbolic approximations of $(\TT^d,\|\cdot\|_{\infty})$ with layers $(L_i)_{i\in\ZZ}$ defined in \cite[p. 47]{MT26}, considering $\lambda=1/r=d^2$ . In fact, we will make use of some results from \cite{MT26} in the special case $\lambda=1/r$.

\begin{proposition}\label{pr19}
    There is $C>0$ such that $$d_{\mathrm{BM}}\big(\lipfree{M_d^n},\ell_1^{\abs{M_d^n}-1}\big)\leq C$$ for every $d\geq 6$ and every $n\in\NN$.
\end{proposition}

\begin{proof}
    By \cite[Definition 4.8 and Theorem 4.9]{MT26}, it is enough to show that for every $\varepsilon>0$, if $N=\varepsilon\ZZ^d\cap\TT^d$, then there is $C\geq1$ satisfying that
    \begin{equation}\label{targeteqBN}
        \big|B_{N}(x,d^2\varepsilon)\big|\leq C\big|B_{N}(x,d^2\varepsilon-3\varepsilon)\big|<\infty\;\;\;\;\;\forall x\in \TT^d
    \end{equation}
    where $B_N(x,R)$ stands for the closed ball in $(N,\|\cdot\|_\infty)$ centered in $x$ with radius $R$.
    Let us show that equation \eqref{targeteqBN} actually holds with $C=e^{7/11}$.
    If $(d^2-3)\varepsilon\geq1/2$ then $$B_N(x,d^2\varepsilon)=B_N(x,(d^2-3)\varepsilon)=N$$ and thus equation \eqref{targeteqBN} holds for $C=1$. Otherwise, $|B_N(x,(d^2-3)\varepsilon)|\geq (2d^2-6)^d$ and $|B_N(x,d^2\varepsilon)|\leq (2d^2+1)^d$. Therefore,
    \begin{align*}
    \frac{\abs{B_N(x,d^2\varepsilon)}}{\abs{B_N(x,d^2\varepsilon-3\varepsilon)}} &\leq \pare{\frac{2d^2+1}{2d^2-6}}^d = \pare{1+\frac{7}{2d^2-6}}^{\frac{2d^2-6}{7}\frac{7d}{2d^2-6}} \\
    &\leq \exp\pare{\frac{7d}{2d^2-6}} \leq \exp\pare{\frac{7}{2d-1}} \leq e^{7/11}
    \end{align*}
    where we use $(1+1/x)^x\leq e$ for $x>0$. 
\end{proof}

We now focus on the Nagata dimension of $M_d$. We analyze it by showing that the space considered in Section \ref{ex1} is biLipschitzly embeddable into $M_d$. Let then $(\ZZ^d,\rho)$ be the space considered in Section \ref{ex1}.

\begin{proposition}\label{pr20}
    There is a universal constant $D>0$ such that, for every $d,N\in\NN$ with $d\geq 6$, the space $([N]^d,\rho)$ is $D$-biLipschitzly embeddable into $M_d$.
\end{proposition}
\begin{proof}
    We first consider $i>0$ such that $d^{2i}>d^{4N}$. The embedding $f:([N]^d,\rho)\to(L_i,\rho_d)$ is then given by $f(x)=\frac{x}{d^{2i}}$. Since $\rho(x,y)\geq 1$ and $\rho_d(f(x),f(y))\geq1$ for every pair of distinct $x,y\in [N]^d$, it will be enough to show that
    \begin{equation}\label{targetrho}
        2\rho(x,y)-6\leq\rho_d(f(x),f(y))\leq 2\rho(x,y)+3.
    \end{equation}
    We denote $v_0^x=\frac{x}{d^{2i}}=f(x)$ and  $v_0^y=\frac{y}{d^{2i}}=f(y)$.

    \medskip

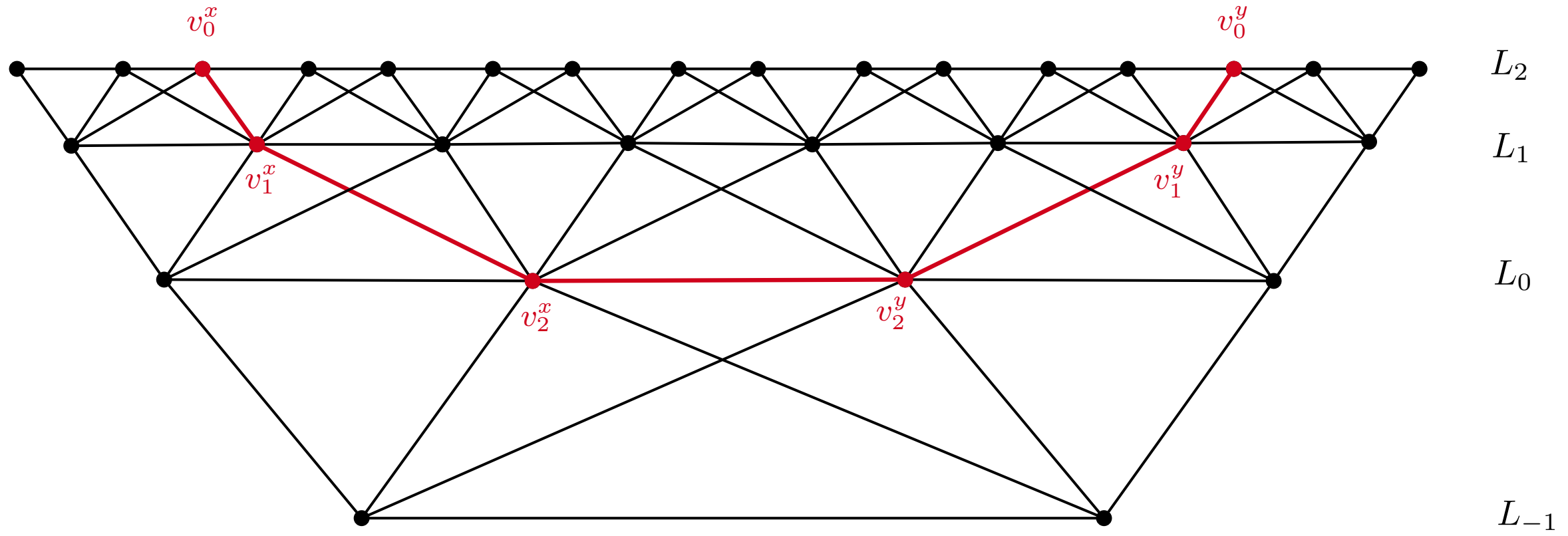
\begin{figure}
    \centering

\tikzset{every picture/.style={line width=0.75pt}} 

\begin{tikzpicture}[x=0.75pt,y=0.75pt,yscale=-1,xscale=1]

\draw  [fill={rgb, 255:red, 0; green, 0; blue, 0 }  ,fill opacity=1 ] (72.67,28.67) .. controls (72.67,27.29) and (73.79,26.17) .. (75.17,26.17) .. controls (76.55,26.17) and (77.67,27.29) .. (77.67,28.67) .. controls (77.67,30.05) and (76.55,31.17) .. (75.17,31.17) .. controls (73.79,31.17) and (72.67,30.05) .. (72.67,28.67) -- cycle ;
\draw  [fill={rgb, 255:red, 0; green, 0; blue, 0 }  ,fill opacity=1 ] (102.67,28.67) .. controls (102.67,27.29) and (103.79,26.17) .. (105.17,26.17) .. controls (106.55,26.17) and (107.67,27.29) .. (107.67,28.67) .. controls (107.67,30.05) and (106.55,31.17) .. (105.17,31.17) .. controls (103.79,31.17) and (102.67,30.05) .. (102.67,28.67) -- cycle ;
\draw    (77.67,28.67) -- (102.67,28.67) ;
\draw  [fill={rgb, 255:red, 0; green, 0; blue, 0 }  ,fill opacity=1 ] (142.67,28.67) .. controls (142.67,27.29) and (143.79,26.17) .. (145.17,26.17) .. controls (146.55,26.17) and (147.67,27.29) .. (147.67,28.67) .. controls (147.67,30.05) and (146.55,31.17) .. (145.17,31.17) .. controls (143.79,31.17) and (142.67,30.05) .. (142.67,28.67) -- cycle ;
\draw    (147.67,28.67) -- (172.67,28.67) ;
\draw    (107.67,28.67) -- (142.67,28.67) ;
\draw  [fill={rgb, 255:red, 0; green, 0; blue, 0 }  ,fill opacity=1 ] (212.17,28.67) .. controls (212.17,27.29) and (213.29,26.17) .. (214.67,26.17) .. controls (216.05,26.17) and (217.17,27.29) .. (217.17,28.67) .. controls (217.17,30.05) and (216.05,31.17) .. (214.67,31.17) .. controls (213.29,31.17) and (212.17,30.05) .. (212.17,28.67) -- cycle ;
\draw  [fill={rgb, 255:red, 0; green, 0; blue, 0 }  ,fill opacity=1 ] (242.17,28.67) .. controls (242.17,27.29) and (243.29,26.17) .. (244.67,26.17) .. controls (246.05,26.17) and (247.17,27.29) .. (247.17,28.67) .. controls (247.17,30.05) and (246.05,31.17) .. (244.67,31.17) .. controls (243.29,31.17) and (242.17,30.05) .. (242.17,28.67) -- cycle ;
\draw    (217.17,28.67) -- (242.17,28.67) ;
\draw    (177.17,28.67) -- (212.17,28.67) ;
\draw  [fill={rgb, 255:red, 0; green, 0; blue, 0 }  ,fill opacity=1 ] (282.17,28.67) .. controls (282.17,27.29) and (283.29,26.17) .. (284.67,26.17) .. controls (286.05,26.17) and (287.17,27.29) .. (287.17,28.67) .. controls (287.17,30.05) and (286.05,31.17) .. (284.67,31.17) .. controls (283.29,31.17) and (282.17,30.05) .. (282.17,28.67) -- cycle ;
\draw  [fill={rgb, 255:red, 0; green, 0; blue, 0 }  ,fill opacity=1 ] (312.17,28.67) .. controls (312.17,27.29) and (313.29,26.17) .. (314.67,26.17) .. controls (316.05,26.17) and (317.17,27.29) .. (317.17,28.67) .. controls (317.17,30.05) and (316.05,31.17) .. (314.67,31.17) .. controls (313.29,31.17) and (312.17,30.05) .. (312.17,28.67) -- cycle ;
\draw    (287.17,28.67) -- (312.17,28.67) ;
\draw    (247.17,28.67) -- (282.17,28.67) ;
\draw  [fill={rgb, 255:red, 0; green, 0; blue, 0 }  ,fill opacity=1 ] (352.17,28.67) .. controls (352.17,27.29) and (353.29,26.17) .. (354.67,26.17) .. controls (356.05,26.17) and (357.17,27.29) .. (357.17,28.67) .. controls (357.17,30.05) and (356.05,31.17) .. (354.67,31.17) .. controls (353.29,31.17) and (352.17,30.05) .. (352.17,28.67) -- cycle ;
\draw  [fill={rgb, 255:red, 0; green, 0; blue, 0 }  ,fill opacity=1 ] (382.17,28.67) .. controls (382.17,27.29) and (383.29,26.17) .. (384.67,26.17) .. controls (386.05,26.17) and (387.17,27.29) .. (387.17,28.67) .. controls (387.17,30.05) and (386.05,31.17) .. (384.67,31.17) .. controls (383.29,31.17) and (382.17,30.05) .. (382.17,28.67) -- cycle ;
\draw    (357.17,28.67) -- (382.17,28.67) ;
\draw    (317.17,28.67) -- (352.17,28.67) ;
\draw  [fill={rgb, 255:red, 0; green, 0; blue, 0 }  ,fill opacity=1 ] (421.67,28.67) .. controls (421.67,27.29) and (422.79,26.17) .. (424.17,26.17) .. controls (425.55,26.17) and (426.67,27.29) .. (426.67,28.67) .. controls (426.67,30.05) and (425.55,31.17) .. (424.17,31.17) .. controls (422.79,31.17) and (421.67,30.05) .. (421.67,28.67) -- cycle ;
\draw  [fill={rgb, 255:red, 0; green, 0; blue, 0 }  ,fill opacity=1 ] (451.67,28.67) .. controls (451.67,27.29) and (452.79,26.17) .. (454.17,26.17) .. controls (455.55,26.17) and (456.67,27.29) .. (456.67,28.67) .. controls (456.67,30.05) and (455.55,31.17) .. (454.17,31.17) .. controls (452.79,31.17) and (451.67,30.05) .. (451.67,28.67) -- cycle ;
\draw    (426.67,28.67) -- (451.67,28.67) ;
\draw    (386.67,28.67) -- (421.67,28.67) ;
\draw  [fill={rgb, 255:red, 0; green, 0; blue, 0 }  ,fill opacity=1 ] (491.67,28.67) .. controls (491.67,27.29) and (492.79,26.17) .. (494.17,26.17) .. controls (495.55,26.17) and (496.67,27.29) .. (496.67,28.67) .. controls (496.67,30.05) and (495.55,31.17) .. (494.17,31.17) .. controls (492.79,31.17) and (491.67,30.05) .. (491.67,28.67) -- cycle ;
\draw  [fill={rgb, 255:red, 0; green, 0; blue, 0 }  ,fill opacity=1 ] (521.67,28.67) .. controls (521.67,27.29) and (522.79,26.17) .. (524.17,26.17) .. controls (525.55,26.17) and (526.67,27.29) .. (526.67,28.67) .. controls (526.67,30.05) and (525.55,31.17) .. (524.17,31.17) .. controls (522.79,31.17) and (521.67,30.05) .. (521.67,28.67) -- cycle ;
\draw    (496.67,28.67) -- (521.67,28.67) ;
\draw    (456.67,28.67) -- (491.67,28.67) ;
\draw    (37.67,28.67) -- (72.67,28.67) ;
\draw    (526.67,28.67) -- (561.67,28.67) ;
\draw  [fill={rgb, 255:red, 0; green, 0; blue, 0 }  ,fill opacity=1 ] (32.67,28.67) .. controls (32.67,27.29) and (33.79,26.17) .. (35.17,26.17) .. controls (36.55,26.17) and (37.67,27.29) .. (37.67,28.67) .. controls (37.67,30.05) and (36.55,31.17) .. (35.17,31.17) .. controls (33.79,31.17) and (32.67,30.05) .. (32.67,28.67) -- cycle ;
\draw  [fill={rgb, 255:red, 0; green, 0; blue, 0 }  ,fill opacity=1 ] (561.67,28.67) .. controls (561.67,27.29) and (562.79,26.17) .. (564.17,26.17) .. controls (565.55,26.17) and (566.67,27.29) .. (566.67,28.67) .. controls (566.67,30.05) and (565.55,31.17) .. (564.17,31.17) .. controls (562.79,31.17) and (561.67,30.05) .. (561.67,28.67) -- cycle ;
\draw  [fill={rgb, 255:red, 0; green, 0; blue, 0 }  ,fill opacity=1 ] (53.17,57.67) .. controls (53.17,56.29) and (54.29,55.17) .. (55.67,55.17) .. controls (57.05,55.17) and (58.17,56.29) .. (58.17,57.67) .. controls (58.17,59.05) and (57.05,60.17) .. (55.67,60.17) .. controls (54.29,60.17) and (53.17,59.05) .. (53.17,57.67) -- cycle ;
\draw  [color={rgb, 255:red, 0; green, 0; blue, 0 }  ,draw opacity=1 ][fill={rgb, 255:red, 0; green, 0; blue, 0 }  ,fill opacity=1 ] (123.17,57.17) .. controls (123.17,55.79) and (124.29,54.67) .. (125.67,54.67) .. controls (127.05,54.67) and (128.17,55.79) .. (128.17,57.17) .. controls (128.17,58.55) and (127.05,59.67) .. (125.67,59.67) .. controls (124.29,59.67) and (123.17,58.55) .. (123.17,57.17) -- cycle ;
\draw  [fill={rgb, 255:red, 0; green, 0; blue, 0 }  ,fill opacity=1 ] (193.17,57.17) .. controls (193.17,55.79) and (194.29,54.67) .. (195.67,54.67) .. controls (197.05,54.67) and (198.17,55.79) .. (198.17,57.17) .. controls (198.17,58.55) and (197.05,59.67) .. (195.67,59.67) .. controls (194.29,59.67) and (193.17,58.55) .. (193.17,57.17) -- cycle ;
\draw  [fill={rgb, 255:red, 0; green, 0; blue, 0 }  ,fill opacity=1 ] (263.17,56.67) .. controls (263.17,55.29) and (264.29,54.17) .. (265.67,54.17) .. controls (267.05,54.17) and (268.17,55.29) .. (268.17,56.67) .. controls (268.17,58.05) and (267.05,59.17) .. (265.67,59.17) .. controls (264.29,59.17) and (263.17,58.05) .. (263.17,56.67) -- cycle ;
\draw  [fill={rgb, 255:red, 0; green, 0; blue, 0 }  ,fill opacity=1 ] (332.67,57.17) .. controls (332.67,55.79) and (333.79,54.67) .. (335.17,54.67) .. controls (336.55,54.67) and (337.67,55.79) .. (337.67,57.17) .. controls (337.67,58.55) and (336.55,59.67) .. (335.17,59.67) .. controls (333.79,59.67) and (332.67,58.55) .. (332.67,57.17) -- cycle ;
\draw  [fill={rgb, 255:red, 0; green, 0; blue, 0 }  ,fill opacity=1 ] (402.67,56.67) .. controls (402.67,55.29) and (403.79,54.17) .. (405.17,54.17) .. controls (406.55,54.17) and (407.67,55.29) .. (407.67,56.67) .. controls (407.67,58.05) and (406.55,59.17) .. (405.17,59.17) .. controls (403.79,59.17) and (402.67,58.05) .. (402.67,56.67) -- cycle ;
\draw  [fill={rgb, 255:red, 0; green, 0; blue, 0 }  ,fill opacity=1 ] (472.67,56.67) .. controls (472.67,55.29) and (473.79,54.17) .. (475.17,54.17) .. controls (476.55,54.17) and (477.67,55.29) .. (477.67,56.67) .. controls (477.67,58.05) and (476.55,59.17) .. (475.17,59.17) .. controls (473.79,59.17) and (472.67,58.05) .. (472.67,56.67) -- cycle ;
\draw  [fill={rgb, 255:red, 0; green, 0; blue, 0 }  ,fill opacity=1 ] (542.67,56.17) .. controls (542.67,54.79) and (543.79,53.67) .. (545.17,53.67) .. controls (546.55,53.67) and (547.67,54.79) .. (547.67,56.17) .. controls (547.67,57.55) and (546.55,58.67) .. (545.17,58.67) .. controls (543.79,58.67) and (542.67,57.55) .. (542.67,56.17) -- cycle ;
\draw    (58.17,57.67) -- (123.17,57.17) ;
\draw    (128.17,57.17) -- (193.17,57.17) ;
\draw    (198.17,57.17) -- (263.17,56.67) ;
\draw    (268.17,56.67) -- (332.67,57.17) ;
\draw    (337.67,57.17) -- (402.67,56.67) ;
\draw    (407.67,56.67) -- (472.67,56.67) ;
\draw    (477.67,56.67) -- (542.67,56.17) ;
\draw    (35.17,28.67) -- (55.67,57.67) ;
\draw    (75.17,28.67) -- (55.67,57.67) ;
\draw [color={rgb, 255:red, 208; green, 2; blue, 27 }  ,draw opacity=1 , very thick]   (105.17,28.67) -- (125.67,57.17) ;
\draw    (75.17,28.67) -- (125.67,57.17) ;
\draw    (105.17,28.67) -- (55.67,57.67) ;
\draw    (145.17,28.67) -- (125.67,57.17) ;
\draw    (145.17,28.67) -- (195.67,57.17) ;
\draw    (195.67,57.17) -- (214.67,28.67) ;
\draw    (244.67,28.67) -- (265.67,56.67) ;
\draw    (195.67,57.17) -- (244.67,28.67) ;
\draw    (265.67,56.67) -- (214.67,28.67) ;
\draw    (284.67,28.67) -- (265.67,56.67) ;
\draw    (314.67,28.67) -- (335.17,57.17) ;
\draw    (314.67,28.67) -- (265.67,56.67) ;
\draw    (284.67,28.67) -- (335.17,57.17) ;
\draw    (354.67,28.67) -- (335.17,57.17) ;
\draw    (384.67,28.67) -- (405.17,56.67) ;
\draw    (354.67,28.67) -- (405.17,56.67) ;
\draw    (384.67,28.67) -- (335.17,57.17) ;
\draw    (424.17,28.67) -- (405.17,56.67) ;
\draw    (424.17,28.67) -- (475.17,56.67) ;
\draw    (454.17,28.67) -- (405.17,56.67) ;
\draw    (454.17,28.67) -- (475.17,56.67) ;
\draw [color={rgb, 255:red, 208; green, 2; blue, 27 }  ,draw opacity=1 , very thick]   (494.17,28.67) -- (475.17,56.67) ;
\draw    (494.17,28.67) -- (545.17,56.17) ;
\draw    (524.17,28.67) -- (475.17,56.67) ;
\draw    (524.17,28.67) -- (545.17,56.17) ;
\draw    (564.17,28.67) -- (545.17,56.17) ;
\draw  [fill={rgb, 255:red, 0; green, 0; blue, 0 }  ,fill opacity=1 ] (88.17,108.17) .. controls (88.17,106.79) and (89.29,105.67) .. (90.67,105.67) .. controls (92.05,105.67) and (93.17,106.79) .. (93.17,108.17) .. controls (93.17,109.55) and (92.05,110.67) .. (90.67,110.67) .. controls (89.29,110.67) and (88.17,109.55) .. (88.17,108.17) -- cycle ;
\draw  [fill={rgb, 255:red, 0; green, 0; blue, 0 }  ,fill opacity=1 ] (227.17,108.67) .. controls (227.17,107.29) and (228.29,106.17) .. (229.67,106.17) .. controls (231.05,106.17) and (232.17,107.29) .. (232.17,108.67) .. controls (232.17,110.05) and (231.05,111.17) .. (229.67,111.17) .. controls (228.29,111.17) and (227.17,110.05) .. (227.17,108.67) -- cycle ;
\draw  [fill={rgb, 255:red, 0; green, 0; blue, 0 }  ,fill opacity=1 ] (367.67,108.17) .. controls (367.67,106.79) and (368.79,105.67) .. (370.17,105.67) .. controls (371.55,105.67) and (372.67,106.79) .. (372.67,108.17) .. controls (372.67,109.55) and (371.55,110.67) .. (370.17,110.67) .. controls (368.79,110.67) and (367.67,109.55) .. (367.67,108.17) -- cycle ;
\draw  [fill={rgb, 255:red, 0; green, 0; blue, 0 }  ,fill opacity=1 ] (506.67,108.67) .. controls (506.67,107.29) and (507.79,106.17) .. (509.17,106.17) .. controls (510.55,106.17) and (511.67,107.29) .. (511.67,108.67) .. controls (511.67,110.05) and (510.55,111.17) .. (509.17,111.17) .. controls (507.79,111.17) and (506.67,110.05) .. (506.67,108.67) -- cycle ;
\draw    (55.67,57.67) -- (90.67,108.17) ;
\draw    (125.67,57.17) -- (90.67,108.17) ;
\draw [color={rgb, 255:red, 208; green, 2; blue, 27 }  ,draw opacity=1 , very thick]   (125.67,57.17) -- (229.67,108.67) ;
\draw    (90.67,108.17) -- (229.67,108.67) ;
\draw    (90.67,108.17) -- (195.67,57.17) ;
\draw    (229.67,108.67) -- (195.67,57.17) ;
\draw    (370.17,108.17) -- (509.17,108.67) ;
\draw [color={rgb, 255:red, 208; green, 2; blue, 27 }  ,draw opacity=1 , very thick]   (229.67,108.67) -- (370.17,108.17) ;
\draw    (229.67,108.67) -- (265.67,56.67) ;
\draw    (229.67,108.67) -- (335.17,57.17) ;
\draw    (370.17,108.17) -- (265.67,56.67) ;
\draw    (370.17,108.17) -- (335.17,57.17) ;
\draw    (370.17,108.17) -- (405.17,56.67) ;
\draw [color={rgb, 255:red, 208; green, 2; blue, 27 }  ,draw opacity=1 , very thick]   (370.17,108.17) -- (475.17,56.67) ;
\draw    (509.17,108.67) -- (405.17,56.67) ;
\draw    (509.17,108.67) -- (545.17,56.17) ;
\draw    (509.17,108.67) -- (475.17,56.67) ;
\draw  [fill={rgb, 255:red, 0; green, 0; blue, 0 }  ,fill opacity=1 ] (162.67,198.17) .. controls (162.67,196.79) and (163.79,195.67) .. (165.17,195.67) .. controls (166.55,195.67) and (167.67,196.79) .. (167.67,198.17) .. controls (167.67,199.55) and (166.55,200.67) .. (165.17,200.67) .. controls (163.79,200.67) and (162.67,199.55) .. (162.67,198.17) -- cycle ;
\draw  [fill={rgb, 255:red, 0; green, 0; blue, 0 }  ,fill opacity=1 ] (442.67,198.17) .. controls (442.67,196.79) and (443.79,195.67) .. (445.17,195.67) .. controls (446.55,195.67) and (447.67,196.79) .. (447.67,198.17) .. controls (447.67,199.55) and (446.55,200.67) .. (445.17,200.67) .. controls (443.79,200.67) and (442.67,199.55) .. (442.67,198.17) -- cycle ;
\draw    (90.67,108.17) -- (165.17,198.17) ;
\draw    (229.67,108.67) -- (165.17,198.17) ;
\draw    (229.67,108.67) -- (445.17,198.17) ;
\draw    (370.17,108.17) -- (445.17,198.17) ;
\draw    (370.17,108.17) -- (165.17,198.17) ;
\draw    (509.17,108.67) -- (445.17,198.17) ;
\draw    (165.17,198.17) -- (445.17,198.17) ;
\draw [color={rgb, 255:red, 0; green, 0; blue, 0 }  ,draw opacity=1 ]   (125.67,57.17) -- (175.17,28.67) ;
\draw  [color={rgb, 255:red, 0; green, 0; blue, 0 }  ,draw opacity=1 ][fill={rgb, 255:red, 0; green, 0; blue, 0 }  ,fill opacity=1 ] (172.67,28.67) .. controls (172.67,27.29) and (173.79,26.17) .. (175.17,26.17) .. controls (176.55,26.17) and (177.67,27.29) .. (177.67,28.67) .. controls (177.67,30.05) and (176.55,31.17) .. (175.17,31.17) .. controls (173.79,31.17) and (172.67,30.05) .. (172.67,28.67) -- cycle ;
\draw  [color={rgb, 255:red, 0; green, 0; blue, 0 }  ,draw opacity=1 ][fill={rgb, 255:red, 0; green, 0; blue, 0 }  ,fill opacity=1 ] (193.17,57.17) .. controls (193.17,55.79) and (194.29,54.67) .. (195.67,54.67) .. controls (197.05,54.67) and (198.17,55.79) .. (198.17,57.17) .. controls (198.17,58.55) and (197.05,59.67) .. (195.67,59.67) .. controls (194.29,59.67) and (193.17,58.55) .. (193.17,57.17) -- cycle ;
\draw [color={rgb, 255:red, 0; green, 0; blue, 0 }  ,draw opacity=1 ]   (175.17,28.67) -- (195.67,57.17) ;
\draw  [color={rgb, 255:red, 208; green, 2; blue, 27 }  ,draw opacity=1 ][fill={rgb, 255:red, 208; green, 2; blue, 27 }  ,fill opacity=1 ] (102.67,28.67) .. controls (102.67,27.29) and (103.79,26.17) .. (105.17,26.17) .. controls (106.55,26.17) and (107.67,27.29) .. (107.67,28.67) .. controls (107.67,30.05) and (106.55,31.17) .. (105.17,31.17) .. controls (103.79,31.17) and (102.67,30.05) .. (102.67,28.67) -- cycle ;
\draw  [color={rgb, 255:red, 208; green, 2; blue, 27 }  ,draw opacity=1 ][fill={rgb, 255:red, 208; green, 2; blue, 27 }  ,fill opacity=1 ] (491.67,28.67) .. controls (491.67,27.29) and (492.79,26.17) .. (494.17,26.17) .. controls (495.55,26.17) and (496.67,27.29) .. (496.67,28.67) .. controls (496.67,30.05) and (495.55,31.17) .. (494.17,31.17) .. controls (492.79,31.17) and (491.67,30.05) .. (491.67,28.67) -- cycle ;
\draw  [color={rgb, 255:red, 208; green, 2; blue, 27 }  ,draw opacity=1 ][fill={rgb, 255:red, 208; green, 2; blue, 27 }  ,fill opacity=1 ] (123.17,57.17) .. controls (123.17,55.79) and (124.29,54.67) .. (125.67,54.67) .. controls (127.05,54.67) and (128.17,55.79) .. (128.17,57.17) .. controls (128.17,58.55) and (127.05,59.67) .. (125.67,59.67) .. controls (124.29,59.67) and (123.17,58.55) .. (123.17,57.17) -- cycle ;
\draw  [color={rgb, 255:red, 208; green, 2; blue, 27 }  ,draw opacity=1 ][fill={rgb, 255:red, 208; green, 2; blue, 27 }  ,fill opacity=1 ] (227.17,108.67) .. controls (227.17,107.29) and (228.29,106.17) .. (229.67,106.17) .. controls (231.05,106.17) and (232.17,107.29) .. (232.17,108.67) .. controls (232.17,110.05) and (231.05,111.17) .. (229.67,111.17) .. controls (228.29,111.17) and (227.17,110.05) .. (227.17,108.67) -- cycle ;
\draw  [color={rgb, 255:red, 208; green, 2; blue, 27 }  ,draw opacity=1 ][fill={rgb, 255:red, 208; green, 2; blue, 27 }  ,fill opacity=1 ] (367.67,108.17) .. controls (367.67,106.79) and (368.79,105.67) .. (370.17,105.67) .. controls (371.55,105.67) and (372.67,106.79) .. (372.67,108.17) .. controls (372.67,109.55) and (371.55,110.67) .. (370.17,110.67) .. controls (368.79,110.67) and (367.67,109.55) .. (367.67,108.17) -- cycle ;
\draw  [color={rgb, 255:red, 208; green, 2; blue, 27 }  ,draw opacity=1 ][fill={rgb, 255:red, 208; green, 2; blue, 27 }  ,fill opacity=1 ] (472.67,56.67) .. controls (472.67,55.29) and (473.79,54.17) .. (475.17,54.17) .. controls (476.55,54.17) and (477.67,55.29) .. (477.67,56.67) .. controls (477.67,58.05) and (476.55,59.17) .. (475.17,59.17) .. controls (473.79,59.17) and (472.67,58.05) .. (472.67,56.67) -- cycle ;

\draw (591.86,190.07) node [anchor=north west][inner sep=0.75pt]    {$L_{-1}$};
\draw (590.53,99.4) node [anchor=north west][inner sep=0.75pt]    {$L_{0}$};
\draw (589.86,51.4) node [anchor=north west][inner sep=0.75pt]    {$L_{1}$};
\draw (589.19,20.07) node [anchor=north west][inner sep=0.75pt]    {$L_{2}$};
\draw (97.67,4.07) node [anchor=north west][inner sep=0.75pt]  [font=\small,color={rgb, 255:red, 208; green, 2; blue, 27 }  ,opacity=1 ]  {$v_{0}^{x}$};
\draw (486.33,3.4) node [anchor=north west][inner sep=0.75pt]  [font=\small,color={rgb, 255:red, 208; green, 2; blue, 27 }  ,opacity=1 ]  {$v_{0}^{y}$};
\draw (119.67,63.57) node [anchor=north west][inner sep=0.75pt]  [font=\small,color={rgb, 255:red, 208; green, 2; blue, 27 }  ,opacity=1 ]  {$v_{1}^{x}$};
\draw (223.67,115.57) node [anchor=north west][inner sep=0.75pt]  [font=\small,color={rgb, 255:red, 208; green, 2; blue, 27 }  ,opacity=1 ]  {$v_{2}^{x}$};
\draw (462.33,63.4) node [anchor=north west][inner sep=0.75pt]  [font=\small,color={rgb, 255:red, 208; green, 2; blue, 27 }  ,opacity=1 ]  {$v_{1}^{y}$};
\draw (357.67,113.07) node [anchor=north west][inner sep=0.75pt]  [font=\small,color={rgb, 255:red, 208; green, 2; blue, 27 }  ,opacity=1 ]  {$v_{2}^{y}$};

\end{tikzpicture}
    \caption{The space $M_d^n$, with a geodesic between $x=v_0^x$ and $y=v_0^y$ highlighted.}
    \label{fig:geodesic}
\end{figure}

    \textbf{First inequality of \eqref{targetrho}:}  By \cite[Lemma 6.2.5]{BS07}, there is a geodesic in $M_d$ connecting $v_0^x$ with $v_0^y$ of the form $\gamma=(v_0^x,v_1^x,\dots,v_m^x,v_m^y,v_{m-1}^y,\dots,v_0^y)$ where $v_j^x,v_j^y\in L_{i-j}$ for every $j=0,\dots,m$, where we potentially incur into an abuse of notation as we allow $v_m^x$ and $v_m^y$ to be the same point. See Figure \ref{fig:geodesic}. Since $\rho_d(v_j^x,v_{j+1}^x)=\rho_d(v_j^y,v_{j+1}^y)=1$ for every $j=0,\dots,m-1$, we have that $v_j^x\in B_\infty(v_{j+1}^x,(\frac{1}{d^2})^{i-j-2})$ and $v_j^y\in B_\infty(v_{j+1}^y,(\frac{1}{d^2})^{i-j-2})$ so that
    \begin{equation}\label{ineqnorms}
        \|v_j^x-v_{j+1}^x\|_\infty\leq (d^2)^{j+2-i}\;\;\text{ , }\;\;\|v_j^y-v_{j+1}^y\|_\infty\leq (d^2)^{j+2-i}\;\;\text{ , }\;\;\|v_m^x-v_m^y\|_\infty\leq 2(d^2)^{m+1-i}.
    \end{equation}
    Therefore, by the triangle inequality,
    \begin{equation}\label{maineqgeo}\begin{aligned}
        \|v_0^x-v_0^y\|_\infty\leq&\|v_m^x-v_m^y\|+ \sum_{j=0}^{m-1}\|v_j^x-v_{j+1}^x\|_\infty+\sum_{j=0}^{m-1}\|v_j^y-v_{j+1}^y\|_\infty\\\stackrel{\eqref{ineqnorms}}{\leq}& 2(d^2)^{m+1-i}+2\sum_{j=0}^{m-1}(d^2)^{j+2-i} \leq (d^2)^{m+2-i}\pare{\frac{2}{d}+\frac{2}{d^2-1}}\leq (d^2)^{m+2-i} .
    \end{aligned}\end{equation}
    Then,
    $$(d^2)^{\rho(x,y)-(i+1)}=\|v_0^x-v_0^y\|_\infty\stackrel{\eqref{maineqgeo}}{\leq}(d^2)^{m+2-i}.$$
    This implies that $\rho(x,y)\leq m+3$ and we are thus done since it is clear from the definition of the geodesic that $$\rho_d(f(x),f(y))=\rho_d(v_0^x,v_0^y)\geq2m\geq 2\rho(x,y)-6.$$

    \medskip
    
    \textbf{Second inequality of \eqref{targetrho}:} We consider, for every $j\in\NN$, elements $v_j^x, v_j^y\in L_{i-j}$ such that $\|v_j^x-v_0^x\|_\infty,\|v_j^y-v_0^y\|_\infty\leq(\frac{1}{d^2})^{i-j}$ (they will be different from the similarly named elements in the previous argument). We claim that $B_\infty(v_j^x,(\frac{1}{d^2})^{i-j-1})\subset B_{\infty}(v_{j+1}^x,(\frac{1}{d^2})^{i-j-2})$ and $B_\infty(v_j^y,(\frac{1}{d^2})^{i-j-1})\subset B_{\infty}(v_{j+1}^y,(\frac{1}{d^2})^{i-j-2})$ so that $\rho_d(v_j^x,v_{j+1}^x)=\rho_d(v_j^y,v_{j+1}^y)=1$ for every $j\in\NN$. Indeed, if $z\in B_\infty(v_j^x,(\frac{1}{d^2})^{i-j-1})$ then
    \begin{align*}
    \norm{z-v_{j+1}^x}_\infty &\leq \norm{z-v_j^x}_\infty + \norm{v_j^x-v_0^x}_\infty + \norm{v_0^x-v_{j+1}^x}_\infty \\
    &\leq d^{2(j+1-i)} + d^{2(j-i)} + d^{2(j+1-i)} \leq d^{2(j+2-i)},
    \end{align*}
    and analogously for $v_j^y,v_{j+1}^y$.

    We claim now that for every $j\geq\rho(x,y)$, we have that $$B_\infty\pare{v_j^x,\pare{\frac{1}{d^2}}^{i-j-1}}\cap B_\infty\pare{v_j^y,\pare{\frac{1}{d^2}}^{i-j-1}}\neq\emptyset$$ and thus $\rho_d(v_j^x,v_j^y)\leq1$. Indeed,
    $$\begin{aligned}\|v_j^x-v_j^y\|_\infty\leq&\|v_j^x-v_0^x\|_\infty+\|v_j^y-v_0^y\|_\infty+\|v_0^x-v_0^y\|_\infty\\\leq&(d^2)^{j-i}+(d^2)^{j-i}+(d^2)^{\rho(x,y)-(i+1)}\stackrel{\rho(x,y)\leq j}{\leq}2(d^2)^{j+1-i}=2\bigg(\frac{1}{d^2}\bigg)^{i-j-1}. \end{aligned}$$
    Therefore, taking $m=\lceil{\rho(x,y)}\rceil<\rho(x,y)+1$,
    $$\rho_d(f(x),f(y))=\rho_d(v_0^x,v_0^y)\leq 2m+1\leq2\rho(x,y)+3.$$
\end{proof}

\begin{proof}[Proof of Theorem \ref{mainthm:lipfree}]
For every $d\geq 6$. we have $\dim_N(M_d)\geq d$ by Propositions \ref{pr14} and \ref{pr20}. Fix a base point $0\in L_0$ for $M_d$, and note that the closed ball $B_{M_d}(0,n)$ is contained in $M_d^n$ for all $n\in\NN$. Thus, by Kalton's decomposition theorem \cite[Proposition 4.3]{Kal04}, there is a universal constant $C_1$ such that $\lipfree{M_d}$ is $C_1$-isomorphic to a $C_1$-complemented subspace of the $\ell_1$-sum $\pare{\bigoplus_{n\in\NN}\lipfree{M_d^{2^n}}}_{\ell_1}$ (see \cite[Lemma 1.2]{AM26} and its proof for the formulation of Kalton's theorem needed here). By Proposition \ref{pr19}, this $\ell_1$-sum is $C_2$-isomorphic to $\pare{\bigoplus_{n\in\NN}\ell_1^{\abs{M_d^{2^n}}-1}}_{\ell_1}=\ell_1$ for a universal constant $C_2$. The same then follows for $\lipfree{M_d}$ by Pe\l czy\'nski's classical results.

Let $M$ be the metric sum of the pointed metric spaces $(M_d,0)$, $d\geq 6$. Then $\dim_N(M)=\infty$ because $M$ contains isometric copies of all $M_d$ and, by \cite[Proposition 3.9]{Wea18}, $\lipfree{M}$ is isometric to $\pare{\bigoplus_{d\geq 6}\lipfree{M_d}}_{\ell_1}$, which is $C$-isomorphic to $\pare{\bigoplus_d\ell_1}_{\ell_1}=\ell_1$. This proves the first statement of Theorem \ref{mainthm:lipfree}.

For the second statement, let $S_n=\bigcup_{d\leq n} M_d^n$, where we identify $M_d$ with its copy in $M$, and similarly for the subset $M_d^n$ of $M_d$. It is clear that $(S_n)$ is increasing and its union is $M$. Note that $S_n$ is isometric to the metric sum of the sets $M_d^n$, $d\leq n$, so the statement follows from Proposition \ref{pr19} and \cite[Proposition 3.9]{Wea18}.
\end{proof}

\section*{Acknowledgments}

Parts of this research were carried out during visits of the first author to the Universidad P\'ublica de Navarra and of the second author to the Universitat Polit\`ecnica de Val\`encia in 2025. Both authors would like to thank the respective hosting institutions for their excellent conditions.

R. J. Aliaga received support from Grants PID2021-122126NB-C33 and PID2025-168101NB-I00 and R. Medina received support from Grants PID2021-122126NB-C31 and PID2025-167660NB-I00, all funded by MICIU/AEI/10.13039/501100011033 and by ERDF/EU.

\section*{Statement on AI use}

The authors acknowledge the use of ChatGPT 5.6-Sol by OpenAI to generate suggestions for keywords and 2020 MSC areas. They take full responsability for the final selection.


\bigskip
\bibliography{bibtfg}
\bibliographystyle{alpha}

\end{document}